\documentclass[10pt,a4paper]{amsart}
\usepackage[%
   pdftitle={Convolution sums of divisor functions},
   pdfauthor={Emmanuel Royer (Universite Montpellier II, I3M, UMR
     CNRS 5149)},
   pdfsubject={quasimodular forms and convolution sums of divisor functions},
   pdfkeywords={quasimodular forms, divisor functions, arithmetical identities},
   colorlinks=false,
   linkcolor=blue,
   citecolor=blue,
   urlcolor=black,
]{hyperref}
\usepackage[T1]{fontenc}
\usepackage{mathrsfs}
\usepackage{amsmath,amsthm,amssymb}
\usepackage{pxfonts}
\theoremstyle{plain}

\theoremstyle{remark}

\author{Emmanuel Royer}
\address{Universit\'e Paul Val\'ery -- Montpellier III\\ MIAp\\
  F--34199 Montpellier cedex 5\\ France}
 \email{emmanuel.royer@m4x.org}
\curraddr{Universit\'e Montpellier II\\ I3M -- UMR CNRS 5149, cc051\\
  F--34095 Montpellier cedex 5\\ France}
 \email{royer@math.univ-montp2.fr}
\title[Convolution sums of divisor functions]{Evaluating convolution sums of the divisor function by quasimodular forms}
\newcommand{\C}{\mathbb{C}}
\def\dd{\mathop{\mathrm{d}\null}%
       \mskip-\thinmuskip\mathord{\null}}
\newcommand{\N}{\mathbb{N}}
\newcommand{\m}[1]{M_{#1}}
\newcommand{\qm}[2]{\widetilde{M}_{#1}^{\leq #2}}
\newcommand{\sldz}{\mathrm{SL}(2,\Z)}
\newcommand{\pk}{\mathscr{H}}
\newcommand{\Q}{\mathbb{Q}}
\def\vec#1{\ensuremath{\mathchoice
           {\mbox{\boldmath$\displaystyle\mathbf{#1}$}}
           {\mbox{\boldmath$\textstyle\mathbf{#1}$}}
           {\mbox{\boldmath$\scriptstyle\mathbf{#1}$}}
           {\mbox{\boldmath$\scriptscriptstyle\mathbf{#1}$}}}}
\newcommand{\Z}{\mathbb{Z}}
\let\plainIm = \Im
\def\Im{\mathop{\plainIm\mkern-2mu \mathit{m}}\nolimits}
\DeclareMathOperator{\lcm}{lcm}
\DeclareMathOperator{\tr}{tr}
\makeatletter
\newcommand*\Pmod[1]{%
  \allowbreak
  \mkern 6mu
  ({\operator@font mod}\,\,#1)%
} 
\makeatother

\theoremstyle{plain}
\newtheorem{theorem}{Theorem}
\newtheorem{lemma}[theorem]{Lemma}
\newtheorem{definition}[theorem]{Definition}
\newtheorem{proposition}[theorem]{Proposition}
\theoremstyle{remark}
\newtheorem{remark}{Remark}
\begin{document}%
\begin{abstract}
We provide a systematic method to compute arithmetic sums including some
previously computed by Alaca, Alaca, Besge, Cheng, Glaisher, Huard, Lahiri, Lemire, Melfi, Ou, Ramanujan,
Spearman and Williams. Our method is based on quasimodular forms. This extension of modular forms has been
constructed by Kaneko \& Zagier.
\end{abstract}
\maketitle


\section{Introduction}
\subsection{Results}
Let $\N$ denote the set of natural numbers and $\N^*=\N\setminus\{0\}$. For $n$ and $j$ in $\N^*$ we set
\[
\sigma_j(n)\coloneqq\sum_{d\mid n}d^j
\]
where $d$ runs through the positive divisors of $n$. If $n\notin\N^*$ we set $\sigma_j(n)=0$. 
Following \cite{MR2173379}, for $N\in\N^*$ we define
\[
W_N(n)\coloneqq\sum_{m<n/N}\sigma_1(m)\sigma_1(n-Nm)
\]
where $m$ runs through the positive integers $<n/N$. We call $W_N$ the convolution of level $N$ (of the divisor function). 
We present a method (introduced in \cite{GT0509205}) to compute these sums
using quasimodular forms.
The comparison between the results we obtain and the ones already obtained may
lead to interesting modular identites (see, for example, remark \ref{rem:imi}). 
We insist on the fact that the only goal of this paper is to present a method
and we recapitulate, in table \ref{tab:www} some of the known results. 
We hope however that some of our results are new (see, for example, theorems
\ref{thm:onze}, \ref{thm:treize}, \ref{thm:quatorze} and proposition \ref{prop:1715}).
Whereas the evaluations of $W_N(n)$ for $N\in\{1,2,3,4\}$ given in \cite{MR1956253}
are elementary and the ones of $W_N(n)$ for $N\in\{5,\dotsc,9\}$
are analytic in nature and use ideas of Ramanujan, our evaluations are on algebraic nature. 

\begin{table}[!ht]\label{tab:www}
\centering
\begin{tabular}{|c|c|c|}
\hline
Level $N$ & Who & Where \\
\hline
$1$ & Besge (Liouville), Glaisher, Ramanujan & \cite{Bes62}, \cite{Gla84}, \cite{Ram16} \\
\hline
$2$, $3$, $4$ & Huard, Ou, Spearman \& Williams & \cite{MR1956253} \\
\hline
$5$, $7$ & Lemire \& Williams & \cite{LeWi05} \\
\hline
$6$ & Alaca \& Williams & \cite{AlWi} \\
\hline
$8$ & Williams & \cite{Wi} \\
\hline
$9$ & Williams & \cite{MR2173379} \\
\hline
$12$ & Alaca, Alaca \& Williams & \cite{AAW12} \\
\hline
$16$ & Alaca, Alaca \& Williams & \cite{AAW16} \\
\hline
$18$ & Alaca, Alaca \& Williams & \cite{AAW18} \\
\hline
$24$ & Alaca, Alaca \& Williams & \cite{AAW24} \\
\hline
\end{tabular}
\caption{Some previous computations of $W_N$}
\end{table}

For $N\in\{5,\dotsc,10\}$, we denote by $\Delta_{4,N}$ the unique cuspidal form
spanning the cuspidal subspace of the modular forms of weight $4$ on $\Gamma_0(N)$
with Fourier expansion\footnote{In this paper, ``Fourier expansion'' always means ``Fourier expansion at the cusp $\infty$''.} 
$\Delta_{4,N}(z)=e^{2\pi iz}+O\left(e^{4\pi iz}\right)$.
We define
\[
\Delta_{4,N}(z)
\eqqcolon
\sum_{n=1}^{+\infty}\tau_{4,N}(n)e^{2\pi inz}.
\]
We also write
\[
\Delta(z)\coloneqq e^{2\pi i  z}\prod_{n=1}^{+\infty}\left[1-e^{2\pi inz}\right]^{24}\\
          \eqqcolon \sum_{n=1}^{+\infty}\tau(n)e^{2\pi inz}
\]
for the unique primitive form of weight $12$ on $\sldz$. 
\begin{theorem}\label{thm:un}
Let $n\in\N^*$, then
\[
W_1(n)
=
\frac{5}{12}\sigma_3(n)
-\frac{n}{2}\sigma_1(n)
+\frac{1}{12}\sigma_1(n),
\]

\[
W_2(n)
=
\frac{1}{12}\sigma_3(n)
+\frac{1}{3}\sigma_3\left(\frac{n}{2}\right)
-\frac{1}{8}n\sigma_1(n)
-\frac{1}{4}n\sigma_1\left(\frac{n}{2}\right)
+\frac{1}{24}\sigma_1(n)
+\frac{1}{24}\sigma_1\left(\frac{n}{2}\right),
\]

\[
W_3(n)
=
\frac{1}{24}\sigma_3(n)
+\frac{3}{8}\sigma_3\left(\frac{n}{3}\right)
-\frac{1}{12}n\sigma_1(n)
-\frac{1}{4}n\sigma_1\left(\frac{n}{3}\right)
+\frac{1}{24}\sigma_1(n)
+\frac{1}{24}\sigma_1\left(\frac{n}{3}\right),
\]

\begin{multline*}
W_4(n)
=
\frac{1}{48}\sigma_3(n)
+\frac{1}{16}\sigma_3\left(\frac{n}{2}\right)
+\frac{1}{3}\sigma_3\left(\frac{n}{4}\right)
-\frac{1}{16}n\sigma_1(n)
-\frac{1}{4}n\sigma_1\left(\frac{n}{4}\right)\\
+\frac{1}{24}\sigma_1(n)
+\frac{1}{24}\sigma_1\left(\frac{n}{4}\right),
\end{multline*}

\begin{multline*}
W_5(n)
=
\frac{5}{312}\sigma_3(n)
+\frac{125}{312}\sigma_3\left(\frac{n}{5}\right)
-\frac{1}{20}n\sigma_1(n)
-\frac{1}{4}n\sigma_1\left(\frac{n}{5}\right)
+\frac{1}{24}\sigma_1(n)
+\frac{1}{24}\sigma_1\left(\frac{n}{5}\right)
\\
-\frac{1}{130}\tau_{4,5}(n),
\end{multline*}

\begin{multline*}
W_6(n)
=
\frac{1}{120}\sigma_3(n)
+\frac{1}{30}\sigma_3\left(\frac{n}{2}\right)
+\frac{3}{40}\sigma_3\left(\frac{n}{3}\right)
+\frac{3}{10}\sigma_3\left(\frac{n}{6}\right)
-\frac{1}{24}n\sigma_1(n)
-\frac{1}{4}n\sigma_1\left(\frac{n}{6}\right)
\\
+\frac{1}{24}\sigma_1(n)
+\frac{1}{24}\sigma_1\left(\frac{n}{6}\right)
-\frac{1}{120}\tau_{4,6}(n),
\end{multline*}

\begin{multline*}
W_7(n)
=
\frac{1}{120}\sigma_3(n)
+\frac{49}{120}\sigma_3\left(\frac{n}{7}\right)
-\frac{1}{28}n\sigma_1(n)
-\frac{1}{4}n\sigma_1\left(\frac{n}{7}\right)
+\frac{1}{24}\sigma_1(n)
+\frac{1}{24}\sigma_1\left(\frac{n}{7}\right)
\\
-\frac{1}{70}\tau_{4,7}(n),
\end{multline*}

\begin{multline*}
W_8(n)
=
\frac{1}{192}\sigma_3(n)
+\frac{1}{64}\sigma_3\left(\frac{n}{2}\right)
+\frac{1}{16}\sigma_3\left(\frac{n}{4}\right)
+\frac{1}{3}\sigma_3\left(\frac{n}{8}\right)
-\frac{1}{32}n\sigma_1(n)
-\frac{1}{4}n\sigma_1\left(\frac{n}{8}\right)
\\
+\frac{1}{24}\sigma_1(n)
+\frac{1}{24}\sigma_1\left(\frac{n}{8}\right)
-\frac{1}{64}\tau_{4,8}(n),
\end{multline*}

\begin{multline*}
W_9(n)
=
\frac{1}{216}\sigma_3(n)
+\frac{1}{27}\sigma_3\left(\frac{n}{3}\right)
+\frac{3}{8}\sigma_3\left(\frac{n}{9}\right)
-\frac{1}{36}n\sigma_1(n)
-\frac{1}{4}n\sigma_1\left(\frac{n}{9}\right)
\\
+\frac{1}{24}\sigma_1(n)
+\frac{1}{24}\sigma_1\left(\frac{n}{9}\right)
-\frac{1}{54}\tau_{4,9}(n).
\end{multline*}

\begin{multline*}
W_{10}(n)
=
\frac{1}{312}\sigma_3(n)
+\frac{1}{78}\sigma_3\left(\frac{n}{2}\right)
+\frac{25}{312}\sigma_3\left(\frac{n}{5}\right)
+\frac{25}{78}\sigma_3\left(\frac{n}{10}\right)
-\frac{1}{40}n\sigma_1(n)
-\frac{1}{4}n\sigma_1\left(\frac{n}{10}\right)
\\
+\frac{1}{24}\sigma_1(n)
+\frac{1}{24}\sigma_1\left(\frac{n}{10}\right)
-\frac{1}{120}\tau_{4,10}(n)
-\frac{3}{260}\tau_{4,5}(n)
-\frac{3}{65}\tau_{4,5}\left(\frac{n}{2}\right).
\end{multline*}
\end{theorem}

For dimensional reasons, the forms $\Delta_{4,N}$ are primitive forms for
$N\in\{5,\dotsc,10\}$, meaning that they are eigenvalues of all the Hecke
operators and that their Fourier expansion begins with $e^{2\pi iz}+O(e^{4\pi
  iz})$. It follows that the arithmetic functions $\tau_{4,N}$ are
multiplicative and satisfy the relation \eqref{eq:hecke} (see below).
Following \cite{MR759465}, one obtains
\begin{align*}
\Delta_{4,5}(z)&=\left[\Delta(z)\Delta(5z)\right]^{1/6},\\
\Delta_{4,6}(z)&=\left[\Delta(z)\Delta(2z)\Delta(3z)\Delta(6z)\right]^{1/12},\\
\Delta_{4,8}(z)&=\left[\Delta(2z)\Delta(4z)\right]^{1/6},\\
\Delta_{4,9}(z)&=\left[\Delta(3z)\right]^{1/3},
\end{align*}
whereas $\Delta_{4,7}$ and $\Delta_{4,10}$ \emph{are not} products of the $\eta$ function. 

However, using \textsc{Magma}\cite{MR1484478} (see \cite{Ste04} for the algorithms based on the computation 
of the spectrum of Hecke operators on modular symbols), 
one can compute their Fourier coefficients (see tables \ref{tab47} and \ref{tab410}).

\begin{remark}\label{rem:imi}
The independant computation of $W_7$ by Lemire \& Williams \cite{LeWi05} implies that
\[
\Delta_{4,7}(z)
=
\left[%
\left(\Delta(z)^2\Delta(7z)\right)^{1/3}%
+%
13\left(\Delta(z)\Delta(7z)\right)^{1/2}%
+%
49\left(\Delta(z)\Delta(7z)^2\right)^{1/3}
\right]^{1/3}.
\]
This provide an alternative method to compute the function $\tau_{4,7}$.
It is likely that, following \cite{LeWi05} to evaluate $W_{10}$ we could get an expression of $\Delta_{4,10}$.
\end{remark}

\begin{table}[!ht]
\centering
\begin{tabular}{|c||c|c|c|c|c|c|c|c|c|c|c|}
\hline
$n$             & $1$  & $2$  & $3$   & $4$  & $5$  & $6$  & $7$  & $8$    & $9$    & $10$  & $11$
\\ \hline
$\tau_{4,7}(n)$ & $1$  & $-1$ & $-2$  & $-7$ & $16$ & $2$  & $-7$ & $15$   & $-23$  & $-16$ & $-8$
\\ \hline
\hline
$n$             & $12$ & $13$ & $14$  & $15$ & $16$ & $17$ & $18$ & $19$   & $20$   & $21$  & $22$ 
\\ \hline
$\tau_{4,7}(n)$ & $14$ & $28$ & $7$ & $-32$  & $41$ & $54$ & $23$ & $-110$ & $-112$ & $14$  & $8$ 
\\ \hline
\end{tabular}
\caption{First Fourier coefficients of $\Delta_{4,7}$}\label{tab47}
\end{table}

\begin{table}[!ht]
\centering
\begin{tabular}{|c||c|c|c|c|c|c|c|c|c|c|c|}
\hline
$n$             & $1$  & $2$  & $3$   & $4$  & $5$  & $6$  & $7$  & $8$    & $9$    & $10$  & $11$
\\ \hline
$\tau_{4,10}(n)$ & $1$ & $2$ & $-8$ & $4$ & $5$ & $-16$ & $-4$ & $8$ & $37$ & $10$ & $12 $
\\ \hline
\hline
$n$             &  $12$ & $13$ & $14$  & $15$ & $16$ & $17$ & $18$ & $19$   & $20$   & $21$  & $22$ 
\\ \hline
$\tau_{4,10}(n)$ & $ -32$ & $-58$ & $-8$ & $-40$ & $16$ & $66$ & $74$ & $-100$ & $20$ & $32$ & $24$
\\ \hline
\end{tabular}
\caption{First Fourier coefficients of $\Delta_{4,10}$}\label{tab410}
\end{table}

In each of our previous examples, we did not leave the field of rational numbers. 
This might not happen, since the primitive forms do not necessarily have rational coefficients. 
However, every evaluation will make use of totally real algebraic numbers for coefficients since
the extension of $\Q$ by the Fourier coefficients of a primitive form is finite and totally 
real \cite[Proposition 1.3]{MR0314801}. To illustrate that fact, we shall evaluate the convolution sum of
level $11$ and $13$. The set of primitive modular forms of weight $4$ on $\Gamma_0(11)$ has two elements. 
The coefficients of these two primitive forms are in $\Q(t)$ where $t$ is a
root of $X^2-2X-2$  (see \S~\ref{cas11} for the use 
of a method founded in \cite{MR1221103}). 
Each primitive form is determined by the beginning of its Fourier expansion:
\begin{align*}
\Delta_{4,11,1}(z)
&=
e^{2\pi iz}+(2-t)e^{4\pi iz}+O(e^{6\pi iz})
\\
\Delta_{4,11,2}(z)
&=
e^{2\pi iz}+te^{4\pi iz}+O(e^{6\pi iz}).
\end{align*} 
We denote by $\tau_{4,11,i}$ the multiplicative function given by the Fourier coefficients of $\Delta_{4,11,i}$.
The two primitive forms, and hence their Fourier coefficients, are conjugate
by $t\mapsto 2-t$  (see, for example, \cite{MR1357209} for the general
result and \S~\ref{cas11} for the 
special case needed here). We give in \S~\ref{cas11} a way to compute the
functions $\tau_{4,11,i}$ for $i\in\{1,2\}$.

\begin{table}[!ht]
\centering
\begin{tabular}{|c||c|c|c|c|c|}
\hline
$n$             &   $1$  & $2$       & $3$       & $4$          & $5$      
\\ \hline
$\tau_{4,11,1}(n)$ & $1$  & $-t + 2$ & $4t - 5$ & $-2t - 2$ & $-8t + 9$   
\\ \hline
\hline
$n$               & $6$          & $7$          & $8$            & $9$           & $10$   
\\ \hline
$\tau_{4,11,1}(n)$ & $5t - 18$ & $4t + 6$ & $10t - 16$ & $-8t + 30$ & $-9t + 34$
\\ \hline
\hline
$n$                &  $11$          &  $12$          & $13$          & $14$           & $15$ 
\\ \hline
$\tau_{4,11,1}(n)$  & $-11$ & $-14t - 6$ & $20t + 20$ & $-6t + 4$ & $12t - 109$ 
\\ \hline
\end{tabular}
\caption{First Fourier coefficients of $\Delta_{4,11,1}$ where $t^2-2t-2=0$}
\end{table}

\begin{theorem}\label{thm:onze}
Let $n\in\N^*$. Then
\begin{multline*}
W_{11}(n)
=
\frac{5}{1464}\sigma_3(n)+\frac{605}{1464}\sigma_3\left(\frac{n}{11}\right)
-\frac{1}{44}n\sigma_1(n)-\frac{1}{4}n\sigma_1\left(\frac{n}{11}\right)
\\
+\frac{1}{24}\sigma_1(n)+\frac{1}{24}\sigma_1\left(\frac{n}{11}\right)
-\frac{2t+43}{4026}\tau_{4,11,1}(n)+\frac{2t-47}{4026}\tau_{4,11,2}(n).
\end{multline*}
\end{theorem}
\begin{remark}
We have 
\[
-\frac{2t+43}{4026}\tau_{4,11,1}(n)+\frac{2t-47}{4026}\tau_{4,11,2}(n)
=
\tr_{\Q(t)/\Q}\left[-\frac{2t+43}{4026}\tau_{4,11,1}(n)\right]\in\Q.
\]
\end{remark}

The set of primitive modular forms of weight $4$ on $\Gamma_0(13)$ has three elements. One of them, we note
$\Delta_{4,13,1}$, has Fourier coefficients in $\Q$. The two others, we note $\Delta_{4,13,2}$ and $\Delta_{4,13,3}$,
have Fourier coefficients in $\Q(u)$ where $u$ is a root of $X^2-X-4$. 
Each of these two primitive form is determined by the beginning of its Fourier expansion:
\begin{align*}
\Delta_{4,13,2}(z)
&=
e^{2\pi iz}+(1-u)e^{4\pi iz}+O(e^{6\pi iz})
\\
\Delta_{4,13,3}(z)
&=
e^{2\pi iz}+ue^{4\pi iz}+O(e^{6\pi iz}).
\end{align*} 
We denote by $\tau_{4,13,i}$ the multiplicative function given by the Fourier coefficients of $\Delta_{4,13,i}$.
The two primitive forms $\Delta_{4,13,2}$ and $\Delta_{4,13,3}$, and hence their Fourier coefficients, are conjugate 
by $u\mapsto 1-u$ (see, for example, \cite{MR1357209}). We compute table
\ref{tab:cinq} by use of \textsc{Magma}.

\begin{table}[!ht]
\centering
\begin{tabular}{|c||c|c|c|c|c|}
\hline
$n$             &   $1$  & $2$       & $3$       & $4$          & $5$      
\\ \hline
$\tau_{4,13,2}(n)$ & $1$  & $-u + 1$ & $3u + 1$ & $-u - 3$ & $-u - 1$   
\\ \hline
\hline
$n$               & $6$          & $7$          & $8$            & $9$           & $10$   
\\ \hline
$\tau_{4,13,2}(n)$ & $-u - 11$ & $-11u + 1$ & $11u - 7$ & $15u + 10$ & $u + 3$
\\ \hline
\hline
$n$                &  $11$          &  $12$          & $13$          & $14$           & $15$ 
\\ \hline
$\tau_{4,13,2}(n)$  & $-12u + 46$ & $-13u - 15$ & $-13$ & $-u + 45$ & $-7u - 13$ 
\\ \hline
\end{tabular}
\caption{First Fourier coefficients of $\Delta_{4,13,2}$ where $u^2-u-4=0$}
\label{tab:cinq}
\end{table}

\begin{theorem}\label{thm:treize}
Let $n\in\N^*$. Then
\begin{multline*}
W_{13}(n)
=
\frac{1}{408}\sigma_3(n)+\frac{169}{408}\sigma_3\left(\frac{n}{13}\right)
-\frac{1}{52}n\sigma_1(n)-\frac{1}{4}n\sigma_1\left(\frac{n}{13}\right)
\\
+\frac{1}{24}\sigma_1(n)+\frac{1}{24}\sigma_1\left(\frac{n}{13}\right)
+\frac{u-6}{442}\tau_{4,13,2}(n)-\frac{u+5}{442}\tau_{4,13,3}(n).
\end{multline*}
\end{theorem}
\begin{remark}
We have 
\[
\frac{u-6}{442}\tau_{4,13,2}(n)-\frac{u+5}{442}\tau_{4,13,3}(n)
=
\tr_{\Q(u)/\Q}\left[\frac{u-6}{442}\tau_{4,13,2}(n)\right]\in\Q.
\]
\end{remark}

The set of primitive modular forms of weight $4$ on $\Gamma_0(14)$ has two elements. Both have coefficients in $\Q$ and we can
distinguish them by the beginning of their Fourier expansion:
\begin{align*}
\Delta_{4,14,1}(z)
&=
e^{2\pi iz}+2e^{4\pi iz}+O(e^{6\pi iz})
\\
\Delta_{4,14,2}(z)
&=
e^{2\pi iz}-2e^{4\pi iz}+O(e^{6\pi iz}).
\end{align*} 
We denote by $\tau_{4,14,i}$ the multiplicative function given by the Fourier coefficients of $\Delta_{4,14,i}$ and give
in \S\ref{cas14} a method to compute these coefficients and get tables \ref{tab:quatorze1} and \ref{tab:quatorze2}.

\begin{table}[!ht]
\centering
\begin{tabular}{|c||c|c|c|c|c|c|c|c|c|c|c|}
\hline
$n$             & $1$  & $2$  & $3$   & $4$  & $5$  & $6$  & $7$  & $8$    & $9$    & $10$  & $11$
\\ \hline
$\tau_{4,14,1}(n)$ & $1$  & $2$  & $-2$  & $4$  & $-12$  & $-4$  & $7$  & $8$  & $-23$  & $-24$  & $48$ 
\\ \hline
\hline
$n$             &  $12$ & $13$ & $14$  & $15$ & $16$ & $17$ & $18$ & $19$   & $20$   & $21$  & $22$  
\\ \hline
$\tau_{4,14,1}(n)$  & $-8$  & $56$  & $14$  & $24$  & $16$  & $-114$  & $-46$  & $2$  & $-48$  & $-14$  & $96$

\\ \hline
\end{tabular}
\caption{First Fourier coefficients of $\Delta_{4,14,1}$}
\label{tab:quatorze1}
\end{table}

\begin{table}[!ht]
\centering
\begin{tabular}{|c||c|c|c|c|c|c|c|c|c|c|c|}
\hline
$n$             & $1$  & $2$  & $3$   & $4$  & $5$  & $6$  & $7$  & $8$    & $9$    & $10$  & $11$
\\ \hline
$\tau_{4,14,2}(n)$  & $1$  & $-2$  & $8$ & $4$  & $-14$  & $-16$  & $-7$  & $-8$  & $37$  & $28$  & $-28$
\\ \hline
\hline
$n$             &  $12$ & $13$ & $14$  & $15$ & $16$ & $17$ & $18$ & $19$   & $20$   & $21$  & $22$ 
\\ \hline
$\tau_{4,14,2}(n)$  & $32$  & $18$  & $14$  & $-112$  & $16$  & $74$  & $-74$  & $80$  & $-56$  & $-56$  & $56$
\\ \hline
\end{tabular}
\caption{First Fourier coefficients of $\Delta_{4,14,2}$}
\label{tab:quatorze2}
\end{table}

\begin{theorem}\label{thm:quatorze}
Let $n\in\N^*$. Then
\begin{multline*}
W_{14}(n)
=
\frac{1}{600}\sigma_3(n)+\frac{1}{150}\sigma_3\left(\frac{n}{2}\right)+\frac{49}{600}\sigma_3\left(\frac{n}{7}\right)
+\frac{49}{150}\sigma_3\left(\frac{n}{14}\right)
-\frac{1}{56}n\sigma_1(n)-\frac{1}{4}n\sigma_1\left(\frac{n}{14}\right)
\\+\frac{1}{24}\sigma_1(n)+\frac{1}{24}\sigma_1\left(\frac{n}{14}\right)
-\frac{3}{350}\tau_{4,7}(n)-\frac{6}{175}\tau_{4,7}\left(\frac{n}{2}\right)
-\frac{1}{84}\tau_{4,14,1}(n)-\frac{1}{200}\tau_{4,14,2}(n).
\end{multline*}
\end{theorem}

\begin{remark}
The fact that for each $N\in\{12,16,18,20,24\}$ there exists only one primitive form of weight $4$ over $\Gamma_0(N)$
and at most one of weight $2$ implies that the only modular forms appearing in the evaluation of the corresponding $W_N$
have rational coefficients.
\end{remark}

Our method, with the introduction of Dirichlet characters, also allows 
to recover another result of Williams \cite[Theorem 1.2]{MR2173379} which extended a result of Melfi 
\cite[Theorem 2, (7)]{MR1628855}. This result is theorem \ref{thm:deux}.
For $b\in\N^*$ and $a\in\{0,\dotsc,b-1\}$, we define
\[
S[a,b](n)\coloneqq 
\sum_{\substack{m=0\\ m\equiv a\Pmod{b}}}^n\sigma_1(m)\sigma_1(n-m).
\]
We compute $S[i,3]$ for $i\in\{0,1,2\}$. Our result uses the primitive Dirichlet character $\chi_3$ defined by
\[
\chi_3(n)
\coloneqq
\begin{cases}
0 & \text{ if $n\equiv 0\pmod{3}$}\\
1 & \text{ if $n\equiv 1\pmod{3}$}\\
-1 & \text{ if $n\equiv -1\pmod{3}$}\\
\end{cases}
\]
for all $n\in\N^*$.
\begin{theorem}\label{thm:deux}
Let $n\in\N^*$, then
\begin{multline*}
S[0,3](n)
=
\frac{11}{72}\sigma_3(n)
+\frac{25}{18}\sigma_3\left(\frac{n}{3}\right)
-\frac{9}{8}\sigma_3\left(\frac{n}{9}\right)
-\frac{1}{4}n\sigma_1(n)
-n\sigma_1\left(\frac{n}{3}\right)
+\frac{3}{4}n\sigma_1\left(\frac{n}{9}\right)\\
+\frac{1}{24}\left[1+\delta(3\mid n)\right]\sigma_1(n)
+
\frac{1}{18}\tau_{4,9}(n),
\end{multline*}
\begin{multline*}
S[1,3](n)
=
\frac{19}{144}\sigma_3(n)
+\frac{1}{48}\chi_3(n)\sigma_3(n)
-\frac{25}{36}\sigma_3\left(\frac{n}{3}\right)
+\frac{9}{16}\sigma_3\left(\frac{n}{9}\right)
\\
-\frac{1}{8}n\sigma_1(n)
-\frac{1}{8}\chi_3(n)n\sigma_1(n)
+\frac{1}{2}n\sigma_1\left(\frac{n}{3}\right)
+\frac{3}{8}\chi_3(n)n\sigma_1\left(\frac{n}{3}\right)
-\frac{3}{8}n\sigma_1\left(\frac{n}{9}\right)
\\
+\frac{1}{24}\delta(3\mid n-1)\sigma_1(n)
+\frac{1}{18}\tau_{4,9}(n)
\end{multline*}
and
\begin{multline*}
S[2,3](n)
=
\frac{19}{144}\sigma_3(n)
-\frac{1}{48}\chi_3(n)\sigma_3(n)
-\frac{25}{36}\sigma_3\left(\frac{n}{3}\right)
+\frac{9}{16}\sigma_3\left(\frac{n}{9}\right)
\\
-\frac{1}{8}n\sigma_1(n)
+\frac{1}{8}\chi_3(n)n\sigma_1(n)
+\frac{1}{2}n\sigma_1\left(\frac{n}{3}\right)
-\frac{3}{8}\chi_3(n)n\sigma_1\left(\frac{n}{3}\right)
-\frac{3}{8}n\sigma_1\left(\frac{n}{9}\right)
\\
+\frac{1}{24}\delta(3\mid n-2)\sigma_1(n)
-\frac{1}{9}\tau_{4,9}(n).
\end{multline*}
where $\delta(3\mid n)$ is $1$ if $3$ divides $n$ and $0$ otherwise. 
\end{theorem}

We next consider convolutions of different divisor sums and recover results of Melfi
\cite[Theorem 2, (9), (10)]{MR1628855} completed by Huard, Ou, Spearman \& Williams \cite[Theorem 6]{MR1956253} and
Cheng \& Williams \cite{CW05}.
We shall use the unique
cuspidal form $\Delta_{8,2}$ 
spanning the cuspidal subspace of the modular forms of weight $8$ on $\Gamma_0(2)$
with Fourier expansion $\Delta_{8,2}(z)=e^{2\pi iz}+O\left(e^{4\pi iz}\right)$.
Using \cite{MR759465}, we have
\[
\Delta_{8,2}(z)=[\eta(z)\eta(2z)]^8.
\]
We define
\[
\Delta_{8,2}(z)
\eqqcolon
\sum_{n=1}^{+\infty}\tau_{8,2}(n)e^{2\pi inz}.
\] 
This is again a primitive form, hence the arithmetic function $\tau_{8,2}$ is multiplicative and
satisfies the relation \eqref{eq:hecke} (see below).

\begin{theorem}\label{thm:trois}
Let $n\in\N^*$. Then
\[
\sum_{k=0}^n\sigma_1(k)\sigma_3(n-k)=
\frac{7}{80}\sigma_5(n)
-\frac{1}{8}n\sigma_3(n)
+\frac{1}{24}\sigma_3(n)
-\frac{1}{240}\sigma_1(n),
\]
\[
\sum_{k<n/2}\sigma_1(n-2k)\sigma_3(k)
=
\frac{1}{240}\sigma_5(n)+\frac{1}{12}\sigma_5\left(\frac{n}{2}\right)-\frac{1}{8}n\sigma_3\left(\frac{n}{2}\right)
+\frac{1}{24}\sigma_3\left(\frac{n}{2}\right)-\frac{1}{240}\sigma_1(n),
\]
\[
\sum_{k<n/2}\sigma_1(k)\sigma_3(n-2k)
=
\frac{1}{48}\sigma_5(n)
+\frac{1}{15}\sigma_5\left(\frac{n}{2}\right)
-\frac{1}{16}n\sigma_3(n)
+\frac{1}{24}\sigma_3(n)
-\frac{1}{240}\sigma_1\left(\frac{n}{2}\right).
\]
Moreover,
\[
\sum_{k=0}^n\sigma_1(k)\sigma_5(n-k)=
\frac{5}{126}\sigma_7(n)
-\frac{1}{12}n\sigma_5(n)
+\frac{1}{24}\sigma_5(n)
+\frac{1}{504}\sigma_1(n),
\]
\begin{multline*}
\sum_{k<n/2}\sigma_1(k)\sigma_5(n-2k)
=
\\
\frac{1}{102}\sigma_7(n)+\frac{32}{1071}\sigma_7\left(\frac{n}{2}\right)-\frac{1}{24}n\sigma_5(n)
+\frac{1}{24}\sigma_5(n)
+\frac{1}{504}\sigma_1\left(\frac{n}{2}\right)-\frac{1}{102}\tau_{8,2}(n),
\end{multline*}
and
\begin{multline*}
\sum_{k<n/2}\sigma_1(n-2k)\sigma_5(k)
=
\\
\frac{1}{2142}\sigma_7(n)+\frac{2}{51}\sigma_7\left(\frac{n}{2}\right)-\frac{1}{12}n\sigma_5\left(\frac{n}{2}\right)
+\frac{1}{24}\sigma_5\left(\frac{n}{2}\right)
+\frac{1}{504}\sigma_1(n)-\frac{1}{408}\tau_{8,2}(n).
\end{multline*}

\end{theorem}
In theorem \ref{thm:trois}, the first and fourth identites are due to Ramanujan \cite{Ram16}. The second and third ones are due
to Huard, Ou, Spearman \& Williams \cite[Theorem 6]{MR1956253}. The fifth and sixth ones are due to Cheng \& Williams \cite{CW05}.
Some other identities of the same type may be found in \cite{CW05} and \cite{Ram16}.

Our method also allows to evaluate sums of Lahiri type
\begin{multline}\label{eq:defS}
S[(a_1,\dotsc,a_r),(b_1,\dotsc,b_r),(N_1,\dotsc,N_r)](n)
\coloneqq
\\
\sum_{\substack{(m_1,\dotsc,m_r)\in\N^r\\ m_1+\dotsc+m_r=n}}
m_1^{a_1}\dotsc m_r^{a_r}\sigma_{b_1}\left(\frac{m_1}{N_1}\right)\dotsc \sigma_{b_r}\left(\frac{m_r}{N_r}\right)
\end{multline}
\cite{MR0020591}, \cite{MR0022566}, \cite[\S 3]{MR1956253} where the $a_i$ are nonnegative integers, the $N_i$ are positive integers and the
$b_i$ are odd positive integers. To simplify the notations, we introduce
\[
S[(a_1,\dotsc,a_r),(b_1,\dotsc,b_r)](n)
\coloneqq
S[(a_1,\dotsc,a_r),(b_1,\dotsc,b_r),(1,\dotsc,1)](n).
\]
For example, we prove the following.
\begin{theorem}\label{thm:quatre}
Let $n\in\N^*$. Then
\begin{multline*}
S[(0,1,1),(1,1,1)](n)
=
\frac{1}{288}n^2\sigma_5(n)
-\frac{1}{72}n^3\sigma_3(n)
+\frac{1}{288}n^2\sigma_3(n)
\\
+\frac{1}{96}n^4\sigma_1(n)
-\frac{1}{288}n^3\sigma_1(n).
\end{multline*}
and if
\begin{align*}
A(n) &= -\frac{48}{5}n^2\sigma_9(n)\\
B(n) &= 128n^3\sigma_7(n)\\
C(n) &= -80n^2\sigma_7(n)-600n^4\sigma_5(n)\\
D(n) &= 648n^3\sigma_5(n)+\frac{8208}{7}n^5\sigma_3(n)\\
E(n) &= -144n^2\sigma_5(n)-\frac{11232}{7}n^4\sigma_3(n)-\frac{3456}{5}n^6\sigma_1(n)\\
F(n) &= 576n^3\sigma_3(n)+\frac{5184}{5}n^5\sigma_1(n)\\
G(n) &= -432n^4\sigma_1(n)-48n^2\sigma_3(n)\\
H(n) &= 48n^3\sigma_1(n)\\
I(n) &= \frac{8}{35}n\tau(n)\\
J(n) &= -\frac{8}{35}\tau(n)
\end{align*}
then 
\begin{multline*}
-24^5
S[(0,0,0,1,1),(1,1,1,1,1)](n)
=\\
A(n)+B(n)+C(n)+D(n)+E(n)+F(n)+G(n)+H(n)+I(n)+J(n).
\end{multline*}
\end{theorem}
The first identity of theorem \ref{thm:quatre} is due to Lahiri \cite[(5.9)]{MR0020591}
and an elementary proof had been given by  Huard, Ou, Spearman \& Williams \cite{MR1956253}. The second
identity is due to Lahiri \cite{MR0022566}.

We continue our evaluations by the more complicated sum
$S[(0,1),(1,1),(2,5)].$ The reason why it is more difficult is that the
underlying space of new cuspidal modular forms has dimension $3$. 

The space of newforms of weight $6$ on $\Gamma_0(10)$ has dimension $3$. Let
$\{\Delta_{6,10,i}\}_{1\leq i\leq 3}$ be the unique basis of primitive forms with
\begin{align*}
\Delta_{6,10,1}(z) & = e^{2\pi iz}+4e^{4\pi iz}+6e^{6\pi iz}+O(e^{8\pi iz}), \\
\Delta_{6,10,2}(z) & =e^{2\pi iz}-4e^{4\pi iz}+24e^{6\pi iz}+O(e^{8\pi iz}), \\
\Delta_{6,10,3}(z) & = e^{2\pi iz}-4e^{4\pi iz}-26e^{6\pi iz}+O(e^{8\pi iz}).
\end{align*}
Again, by \cite{MR759465}, we know that these functions \emph{are not} products of the $\eta$ function.
We denote by $\tau_{6,10,i}(n)$ the $n$th Fourier coefficient of
$\Delta_{6,10,i}$. The functions $\tau_{6,10,i}$ are multiplicative and we
show in \S~\ref{pa:prec} how to establish tables \ref{tab:six},
\ref{tab:sept} and \ref{tab:huit}. 
 
\begin{table}[!ht]
\centering
\begin{tabular}{|c||c|c|c|c|c|c|c|c|c|c|c|}
\hline
$n$             & $1$  & $2$  & $3$   & $4$  & $5$  & $6$  & $7$  & $8$    & $9$    & $10$  & $11$
\\ \hline
$\tau_{6,10,1}(n)$ & $1$ & $ 4$ & $ 6$ & $ 16$ & $ -25$ & $ 24$ & $ -118$ & $ 64$ & $ -207$ & $ -100$ & $ 192$
\\ \hline
\hline
$n$             &  $12$ & $13$ & $14$  & $15$ & $16$ & $17$ & $18$ & $19$   & $20$   & $21$  & $22$  
\\ \hline
$\tau_{6,10,1}(n)$ & $96$ & $ 1106$ & $ -472$ & $ -150$ & $ 256$ & $ 762$ & $ -828$ & $ -2740$ & $ -400$ & $ -708$ & $ 768$

\\ \hline
\end{tabular}
\caption{First Fourier coefficients of $\Delta_{6,10,1}$}
\label{tab:six}
\end{table}

\begin{table}[!ht]
\centering
\begin{tabular}{|c||c|c|c|c|c|c|c|c|c|c|c|}
\hline
$n$             & $1$  & $2$  & $3$   & $4$  & $5$  & $6$  & $7$  & $8$    & $9$    & $10$  & $11$
\\ \hline
$\tau_{6,10,2}(n)$ & $1$ & $ -4$ & $ 24$ & $ 16$ & $ 25$ & $ -96$ & $ -172$ & $ -64$ & $ 333$ & $ -100$ & $ 132$
\\ \hline
\hline
$n$             &  $12$ & $13$ & $14$  & $15$ & $16$ & $17$ & $18$ & $19$   & $20$   & $21$  & $22$ 
\\ \hline
$\tau_{6,10,2}(n)$ & $384$ & $ -946$ & $ 688$ & $ 600$ & $ 256$ & $ -222$ & $ -1332$ & $ 500$ & $ 400$ & $ -4128$ & $ -528$
\\ \hline
\end{tabular}
\caption{First Fourier coefficients of $\Delta_{6,10,2}$}
\label{tab:sept}
\end{table}

\begin{table}[!ht]
\centering
\begin{tabular}{|c||c|c|c|c|c|c|c|c|c|c|c|}
\hline
$n$             & $1$  & $2$  & $3$   & $4$  & $5$  & $6$  & $7$  & $8$    & $9$    & $10$  & $11$
\\ \hline
$\tau_{6,10,3}(n)$ & $1$ & $ -4$ & $ -26$ & $ 16$ & $ -25$ & $ 104$ & $ -22$ & $ -64$ & $ 433$ & $ 100$ & $ -768$
\\ \hline
\hline
$n$             &  $12$ & $13$ & $14$  & $15$ & $16$ & $17$ & $18$ & $19$   & $20$   & $21$  & $22$ 
\\ \hline
$\tau_{6,10,3}(n)$ & $-416$ & $ -46$ & $ 88$ & $ 650$ & $ 256$ & $ 378$ & $ -1732$ & $ 1100$ & $ -400$ & $ 572$ & $ 3072$
\\ \hline
\end{tabular}
\caption{First Fourier coefficients of $\Delta_{6,10,3}$}
\label{tab:huit}
\end{table}

We also need the unique primitive form
\[
\Delta_{6,5}(z)\eqqcolon\sum_{n=1}^{+\infty}\tau_{6,5}(n)e^{2\pi inz}
\]
of weight $6$ on $\Gamma_0(5)$. It is not a product of the $\eta$ function,
and we show in \S~\ref{pa:prec} how to establish table~\ref{tab:neuf}.

\begin{table}[!ht]
\centering
\begin{tabular}{|c||c|c|c|c|c|c|c|c|c|c|c|}
\hline
$n$             & $1$  & $2$  & $3$   & $4$  & $5$  & $6$  & $7$  & $8$    & $9$    & $10$  & $11$
\\ \hline
$\tau_{6,5}(n)$ & $ 1$ & $2$ & $-4$ & $-28$ & $25$ & $-8$ & $192$ & $-120$ & $-227$ & $50$ & $-148$
\\ \hline
\hline
$n$             &  $12$ & $13$ & $14$  & $15$ & $16$ & $17$ & $18$ & $19$   & $20$   & $21$  & $22$  
\\ \hline
$\tau_{6,5}(n)$ & $112$ & $286$ & $384$ & $-100$ & $656$ & $-1678$ & $-454$ & $1060$ & $-700$ & $-768$ & $-296 $
\\ \hline
\end{tabular}
\caption{First Fourier coefficients of $\Delta_{6,5}$}
\label{tab:neuf}
\end{table}

\begin{proposition}\label{prop:1721}
Let $n\in\N^*$. Define
\begin{align*}
A(n) &=
\frac{12}{13}n\sigma_3(n)+\frac{48}{13}n\sigma_3\left(\frac{n}{2}\right)+\frac{300}{13}n\sigma_3\left(\frac{n}{5}\right)
+\frac{1200}{13}n\sigma_3\left(\frac{n}{10}\right)\\
B(n) &=
-\frac{48}{5}n^2\sigma_1\left(\frac{n}{2}\right)-48n^2\sigma_1\left(\frac{n}{5}\right)\\
C(n)
&=
24n\sigma_1\left(\frac{n}{5}\right)\\
D(n)
&=
\frac{12}{5}n\tau_{4,10}(n)-\frac{216}{65}n\tau_{4,5}(n)-\frac{864}{65}n\tau_{4,5}\left(\frac{n}{2}\right)\\
E(n)
&=
\frac{108}{35}\tau_{6,5}(n)+\frac{864}{35}\tau_{6,5}\left(\frac{n}{2}\right)
-\frac{24}{5}\tau_{6,10,1}(n)+\frac{12}{7}\tau_{6,10,2}(n).
\end{align*}
Then
\[
5\times 24^2
\sum_{\substack{(a,b)\in\N^{*2}\\ 2a+5b=n}}b\sigma_1(a)\sigma_1(b)
=
A(n)+B(n)+C(n)+D(n)+E(n).
\]
\end{proposition}

We shall now evaluate $S[(1,1),(1,1),(1,5)]$ since it constitutes an exemple leaving the rational field. 
Let $v$ be one of the two roots of $X^2-20X+24$. There exist three primitive forms
of weight $8$ on $\Gamma_0(5)$ determined by the beginning of their Fourier expansion:
\begin{align*}
\Delta_{8,5,1}(z)
&=
e^{2\pi iz}-14e^{4\pi iz}+O(e^{6\pi iz})
\\
\Delta_{8,5,2}(z)
&=
e^{2\pi iz}+(20-v)e^{4\pi iz}+O(e^{6\pi iz})
\\
\Delta_{8,5,3}(z)
&=
e^{2\pi iz}+ve^{4\pi iz}+O(e^{6\pi iz}).
\end{align*}
The function $\Delta_{8,5,3}$ is obtained from $\Delta_{8,5,2}$ by the conjugation ($v\mapsto 20-v$) of $\Q(v)$ 
on the Fourier coefficients. 
We denote by $\tau_{8,5,i}$ the multiplicative function given by the Fourier
coefficients of $\Delta_{8,5,i}$ and show
in \S~\ref{pa:eight} how to compute them.
\begin{table}[!ht]
\centering
\begin{tabular}{|c||c|c|c|c|c|c|c|c|c|}
\hline
$n$             & $1$  & $2$  & $3$   & $4$  & $5$  & $6$  & $7$  & $8$    & $9$    
\\ \hline
$\tau_{8,5,1}(n)$ & $1$ & $-14$ & $-48$ & $68$ & $125$ & $672$ & $-1644$ & $840$ & $117$
\\ \hline
\hline
$n$             & $10$  & $11$ &  $12$ & $13$ & $14$  & $15$ & $16$ & $17$ & $18$ 
\\ \hline
$\tau_{8,5,1}(n)$ & $-1750$ & $172$ & $-3264$ & $3862$ & $23016$ & $-6000$ & $-20464$ & $-12254$ & $-1638$
\\ \hline
\end{tabular}
\caption{First Fourier coefficients of $\Delta_{8,5,1}$}
\end{table}

\begin{table}[!ht]
\centering
\begin{tabular}{|c||c|c|c|c|c|}
\hline
$n$             &   $1$  & $2$       & $3$       & $4$          & $5$      
\\ \hline
$\tau_{8,5,2}(n)$ & $1$  & $-v + 20$ & $8v - 70$ & $-20v + 248$ & $-125$   
\\ \hline
\hline
$n$               & $6$          & $7$          & $8$            & $9$           & $10$   
\\ \hline
$\tau_{8,5,2}(n)$ & $70v - 1208$ & $-56v + 510$ & $-120v + 1920$ & $160v + 1177$ & $125v - 2500$
\\ \hline
\hline
$n$                &  $11$          &  $12$          & $13$          & $14$           & $15$ 
\\ \hline
$\tau_{8,5,2}(n)$  & $-400v + 6272$ & $184v - 13520$ & $608v - 4310$ & $-510v + 8856$ & $-1000v + 8750$ 
\\ \hline
\end{tabular}
\caption{First Fourier coefficients of $\Delta_{8,5,2}$ where $v^2-20v+24=0$}
\end{table}

\begin{proposition}\label{prop:1715}
Let $n\in\N^*$. Define
\begin{align*}
A(n) &= \frac{24}{13}n^2\sigma_3(n)+\frac{600}{13}n^2\sigma_3\left(\frac{n}{5}\right)\\
B(n) &= -\frac{24}{5}n^3\sigma_1(n)-24n^3\sigma_1\left(\frac{n}{5}\right)\\
C(n) &= -\frac{288}{325}n^2\tau_{4,5}(n)\\
D(n) &= \frac{792+12v}{475}\tau_{8,5,2}(n)+\frac{1032-12v}{475}\tau_{8,5,3}(n).
\end{align*}
Then
\[
5\times 24^2
\sum_{\substack{(a,b)\in\N^{*2}\\ a+5b=n}}ab\sigma_1(a)\sigma_1(b)
=
A(n)+B(n)+C(n)+D(n).
\]
\end{proposition}
\begin{remark}
The two terms in the right hand side of the definition of $D(n)$ in proposition \ref{prop:1715} 
being conjugate, we have 
\[
D(n)=\tr_{\Q(v)/\Q}\left[\frac{792+12v}{475}\tau_{8,5,2}(n)\right]\in\Q.
\]
\end{remark}
To stay in the field of rational numbers, we could have used the fundamental fact that, 
for every even $k>0$ and every integer $N\geq 1$, the space of cuspidal forms of weight
$k$ on $\Gamma_0(N)$ has a basis whose elements have a Fourier expansion with integer coefficients 
\cite[Theorem 3.52]{MR1291394}. However, the coefficients of these Fourier expansions are often not multiplicative: this 
is a good reason to leave $\Q$.

\begin{remark}
If $\tau_{\ast}$ is one of our $\tau$ functions, its values are the Fourier coefficients of a primitive form (of weight $k$ on
$\Gamma_0(N)$ say). It therefore satisfies the
following multiplicativity relation
\begin{equation}\label{eq:hecke}
\tau_{\ast}(mn)
=
\sum_{\substack{d\mid (m,n)\\ (d,N)=1}}\mu(d)d^{k-1}
\tau_{\ast}\left(\frac{m}{d}\right)
\tau_{\ast}\left(\frac{n}{d}\right).
\end{equation}
\end{remark}

\begin{remark}
When we found some, we give some expression to compute the various Fourier
coefficients that we need.
This is however somewhat {\em ad hoc} and never needed since we only need to
compute a finite number of coefficients and can then use the algorithms
provided by algorithmic number theory \cite{MR1484478}, \cite{Ste04}.
\end{remark}
\noindent\textbf{Thanks --}
While a preliminary version of this paper was in circulation,
K.S. Williams kindly informed me of the papers
\cite{AlWi}, \cite{LeWi05}, \cite{Wi}, \cite{AAW12}, \cite{AAW16}, \cite{AAW18}, \cite{AAW24} and \cite{CW05}.
I respectfully thank him for having made these papers available
to me. The final version of this paper was written during my
stay at the Centre de Recherches Math\'ematiques de Montr\'eal which provided me
with very good working conditions. I thank Andrew Granville and Chantal David for their invitation.

\subsection{Method}
Since our method is based on quasimodular forms (anticipated by Rankin \cite{MR0082563} and formally introduced by 
Kaneko \& Zagier in \cite{MR1363056}), we briefly recall the basics on these functions, referring to \cite{MR2186573}
and \cite{GT0509205} for the details. 

Define
\[
\Gamma_0(N)=\left\{
\begin{pmatrix} a & b\\ c & d\end{pmatrix} \colon (a,b,c,d)\in\Z^4,\, ad-bc=1,\, N\mid c
\right\}
\]
for all integers $N\geq 1$.  In particular, $\Gamma_0(1)$ is
$\sldz$.  Denote by $\pk$ the Poincar\'e upper half
plane:
\[
\pk=\{ z\in\C \colon \Im z>0\}.
\]
\begin{definition}
Let $N\in\N$, $k\in\N^*$ and $s\in\N^*$. A holomorphic function
\[
f \colon \pk \to \C
\]
is a quasimodular form of weight $k$, depth $s$ on $\Gamma_0(N)$ if
there exist holomorphic functions $f_0$, $f_1$, $\dotsc$, $f_s$ on
$\pk$ such that
\begin{equation}\label{eq:cqm}
(cz+d)^{-k}
f\left(\frac{az+b}{cz+d}\right)
=
\sum_{i=0}^sf_i(z)\left(\frac{c}{cz+d}\right)^i
\end{equation}
for all $\bigl(\begin{smallmatrix} a & b\\ c &
d\end{smallmatrix}\bigr)\in\Gamma_0(N)$ and such that $f_s$ is
holomorphic at the cusps and not identically vanishing.  By
convention, the $0$ function is a quasimodular form of depth $0$ for
each weight.
\end{definition}
Here is what is meant by the requirement for $f_s$ to be holomorphic at
the cusps.  One can show \cite[Lemme 119]{MR2186573} that if $f$ satisfies
the quasimodularity condition \eqref{eq:cqm}, then $f_s$ satisfies the
modularity condition
\[
(cz+d)^{-(k-2s)}
f_s\left(\frac{az+b}{cz+d}\right)
=
f_s(z)
\]
for all $\bigl(\begin{smallmatrix} a & b\\ c &
d\end{smallmatrix}\bigr)\in\Gamma_0(N)$.  Asking $f_s$ to be
holomorphic at the cusps is asking that, for all
$M=\bigl(\begin{smallmatrix} \alpha & \beta\\ \gamma &
\delta\end{smallmatrix}\bigr)\in\Gamma_0(1)$, the function
\[
z\mapsto(\gamma z+\delta)^{-(k-2s)}f_s\left(\frac{\alpha z+\beta}{\gamma z+\delta}\right)
\]
has a Fourier expansion of the form
\[
\sum_{n=0}^{+\infty}\widehat{f}_{s,M}(n)e\left(\frac{nz}{u_M}\right)
\]
where
\[
u_M=\inf\{
u\in\N^* \colon T^u\in M^{-1}\Gamma_0(N)M
\}.
\]
In other words, $f_s$ is automatically a modular function and is
required to be more than that, a modular form of weight $k-2s$
on $\Gamma_0(N)$.  It follows that if $f$ is a quasimodular form of
weight $k$ and depth $s$, non identically vanishing, then $k$ is even
and $s\leq k/2$.

A fundamental quasimodular form is the Eisenstein series of weight $2$ defined by
\[
E_2(z)=1-24\sum_{n=1}^{+\infty}\sigma_1(n)e^{2\pi inz}.
\]
It is a quasimodular form of weight $2$, depth $1$ on $\Gamma_0(1)$ (see, for example, \cite[Chapter 7]{MR0498338}).

We shall denote by $\qm{k}{s}[\Gamma_0(N)]$ the space of quasimodular forms of weight $k$, depth $\leq s$
on $\Gamma_0(N)$ and $\m{k}[\Gamma_0(N)]=\qm{k}{0}[\Gamma_0(N)]$ the space of modular forms of weight $k$
on $\Gamma_0(N)$. The space $\qm{k}{k/2}[\Gamma_0(N)]$ is graded by the depth.

Our method for theorem \ref{thm:un} is to remark that the function
\begin{align*}
H_N(z)
&= E_2(z)E_2(Nz)\\
&=
1-24\sum_{n=1}^{+\infty}\left[\sigma_1(n)+\sigma_1\left(\frac{n}{N}\right)\right]e^{2\pi inz}
+576\sum_{n=1}^{+\infty}W_N(n)e^{2\pi inz}
\end{align*}
is a quasimodular form of weight $4$, depth $2$ on $\Gamma_0(N)$ that we linearise using the 
following lemma.
\begin{lemma}\label{lem:0040}
Let $k\geq 2$ even. Then
\[
\qm{k}{k/2}[\Gamma_{0}(N)]
=
\bigoplus_{i=0}^{k/2-1}D^i\m{k-2i}[\Gamma_0(N)]\oplus \C D^{k/2-1}E_2.
\]
\end{lemma} 
We have set
\[
D\coloneqq \frac{1}{2\pi i }\frac{\dd}{\dd z}.
\]

Let $\{B_k\}_{k\in\N}$ be the sequence of rational numbers defined by its exponential generating function
\[
\frac{t}{e^t-1}=\sum_{k=0}^{+\infty}B_k\frac{t^k}{k!}.
\]
We shall use the Eisenstein series to express the basis we need:
\[
E_{k,N}(z)\coloneqq 1-\frac{2k}{B_k}\sum_{n=1}^{+\infty}\sigma_{k-1}(n)e^{2\pi inNz}
\in
\m{k}[\Gamma_0(N)]
\]
for all $k\in 2\N^*+2$, $N\in\N^*$. If $N=1$ we simplify by writing $E_k\coloneqq E_{k,N}$.
For weight $2$ forms, we shall need
\[
\Phi_{a,b}(z)
=
\frac{1}{b-a}
\left[bE_2(bz)-aE_2(az)\right]
\in\m{2}[\Gamma_0(b)]
\]
for all $b>1$ and $a\mid b$.

Let $\chi$ be a Dirichlet character.  If $f$
satisfies all of what is needed to be a quasimodular form except
\eqref{eq:cqm} being replaced by
\[
(cz+d)^{-k}
f\left(\frac{az+b}{cz+d}\right)
=
\chi(d)\sum_{i=0}^nf_i(z)\left(\frac{c}{cz+d}\right)^i,
\]
then one says that $f$ is a quasimodular form of weight $k$, depth
$s$ and character $\chi$ on $\Gamma_0(N)$. (In particular, we require
$f_s$ to be a modular form of character $\chi$). We denote by
$\qm{k}{s}[\Gamma_0(N),\chi]$ the vector space of quasimodular forms
of weight $k$, depth $\leq s$ and character $\chi$ on $\Gamma_0(N)$.
If $\chi=\chi_0$ is a principal character of modulus dividing $N$,
then $\qm{k}{s}[\Gamma_0(N),\chi]=\qm{k}{s}[\Gamma_0(N)]$.

If $f\in\qm{k}{s}[\Gamma_0(N)]$, then $f$ has a Fourier expansion
with coefficients $\{\widehat{f}(n)\}_{n\in\N}$. We define the
twist of $f$ by the Dirichlet character $\chi$ as
\[
f\otimes\chi(z)
=
\sum_{n=0}^{+\infty}\chi(n)\widehat{f}(n)e^{2\pi inz}. 
\]
In \cite[Proposition 9]{GT0509205}, we proved the following proposition.
\begin{proposition}\label{prop:twist}
Let $\chi$ be a primitive Dirichlet character of conductor $m$.  Let
$f$ be a quasimodular form of weight $k$ and depth $s$ on
$\Gamma_0(N)$.  Then $f\otimes\chi$ is a quasimodular form of weight
$k$, depth less than or equal to $s$ and character $\chi^2$ on
$\Gamma_0\left(\lcm(N,m^2)\right)$.
\end{proposition}
\begin{remark}
The condition of primitivity of the character may be replaced by the condition of non vanishing
of its Gauss sum.
\end{remark}

The proof of theorem \ref{thm:deux} follows from the linearisation of $E_2\cdot E_2\otimes\chi_3$.

Theorems \ref{thm:trois} and \ref{thm:quatre} follow from the linearisation of derivatives of forms of type $E_jE_{k,N}$.

\subsection{Generalisation of the results}
For $N\geq 1$ and $k\geq 2$, let $A^*_{N,k}$ be the set of triples $(\psi,\phi,t)$ such that
$\psi$ is a primitive Dirichlet character of modulus $L$, $\phi$ is a primitive Dirichlet character of modulus $M$
and $t$ is an integer such that $tLM\mid N$ (and $tLM\neq 1$ if $k=2$) with the extra condition
\begin{equation}\label{eq:extra}
\psi\phi(n)
=
\begin{cases}
1 & \text{ if $(n,N)=1$}\\
0 & \text{otherwise}
\end{cases}
\qquad\text{($n\in\mathbb{N^*}$).}
\end{equation} 
We write $\mathbf{1}$ for the primitive character of modulus $1$ (the constant function
$n\mapsto 1$).
We extend the definition of $\sigma_k$: for $k$ and $n$ in $\N^*$ we set
\[
\sigma_k^{\psi,\phi}(n)\coloneqq\sum_{d\mid n}\psi\left(\frac{n}{d}\right)\phi(d)d^k
\]
where $d$ runs through the positive divisors of $n$. If $n\notin\N^*$ we set $\sigma_k^{\psi,\phi}(n)=0.$
If $M$ is the modulus of the primitive character $\phi$, we define the sequence $\{B_k^\phi\}_{k\in\N}$
by its exponential generating function
\[
\sum_{c=0}^{M-1}\phi(c)
\frac{te^{ct}}{e^{Mt}-1}=\sum_{k=0}^{+\infty}B_k^{\phi}\frac{t^k}{k!}.
\]

For any $(\psi,\phi,t)\in A^*_{N,k}$, define
\[
E_{k}^{\psi,\phi}(z)
\coloneqq
\delta(\psi=\mathbf{1})
-\frac{2k}{B_{k}^\phi}\sum_{n=1}^{+\infty}\sigma_{k-1}^{\psi,\phi}(n)e^{2\pi inz}
\]
and
\[
E_{k,t}^{\psi,\phi}(z)
\coloneqq
\begin{cases}
E_{k}^{\psi,\phi}(tz) & \text{if $(k,\psi,\phi)\neq(2,\mathbf{1},\mathbf{1})$}\\
E_{2}^{\mathbf{1},\mathbf{1}}(z)-tE_{2}^{\mathbf{1},\mathbf{1}}(tz) &\text{otherwise}
\end{cases}
\]
where $\delta(\psi=\mathbf{1})$ is $1$ if $\psi=\mathbf{1}$ and $0$ otherwise.

For $N\geq 1$ and $k\geq 2$ even, the set
\[
\left\{
E_{k,t}^{\psi,\phi} \colon (\psi,\phi,t)\in A^*_{N,k}
\right\}
\]
is a basis for the orthogonal subspace (called Eisenstein subspace, the scalar product being the Petersson one) 
of the cuspidal subspace $S_k[\Gamma_0(N)]$ of $\m{k}[\Gamma_0(N)]$ \cite[Chapter 4]{MR2112196}.

Moreover, by Atkin-Lehner-Li theory \cite[Chapter 5]{MR2112196}, a basis for $S_k[\Gamma_0(N)]$ is
\[
\bigcup_{\substack{(d,M)\in\N^* \\ dM\mid N}}
\alpha_d\left(H_k^*[\Gamma_0(M)]\right)
\]
where $\alpha_d$ is
\[
\begin{array}{ccccc}
\alpha_d & \colon & \m{k}[\Gamma_0(M)] & \to & \m{k}[\Gamma_0(M)]\\
         &        &  f                 & \mapsto & [z\mapsto f(dz)]
\end{array}
\]
and $H_k^*[\Gamma_0(M)]$ is the set of primitive forms of weight $k$ on $\Gamma_0(M)$.

A corollary is the following generalisation of theorems \ref{thm:un} and \ref{thm:trois}. If $f$ is a modular form, we denote by $\{\widehat{f}(n)\}_{n\in\N}$
the sequence of its Fourier coefficients.
\begin{proposition}\label{prop:son}
Let $N\geq 1$. There exist scalars $a_{\psi,\phi,t}$, $a_{M,d,f}$ and $a$ such that, for all $n\geq 1$, we have
\begin{multline*}
W_N(n)
=
\sum_{(\psi,\phi,t)\in A_{N,4}^*}a_{\psi,\phi,t}\sigma_3^{\psi,\phi}\left(\frac{n}{t}\right)
+
\sum_{(\psi,\phi,t)\in A_{N,2}^*}a_{\psi,\phi,t}n\sigma_1^{\psi,\phi}\left(\frac{n}{t}\right)
+an\sigma_1(n)
\\
+\sum_{\substack{(d,M)\in\N^* \\ dM\mid N}}
\sum_{f\in H_4^*[\Gamma_0(M)]}a_{M,d,f}\widehat{f}\left(\frac{n}{d}\right)
+\sum_{\substack{(d,M)\in\N^* \\ dM\mid N}}
\sum_{f\in H_2^*[\Gamma_0(M)]}a_{M,d,f}n\widehat{f}\left(\frac{n}{d}\right)
\\
+\frac{1}{24}\sigma_1(n)
+\frac{1}{24}\sigma_1\left(\frac{n}{N}\right).
\end{multline*}
More generally, for any $N\geq 1$ and any even $\ell\geq 4$, the arithmetic functions
\[
n\mapsto \sum_{k<n/N}\sigma_1(n-kN)\sigma_{\ell-1}(k)-\frac{B_\ell}{2\ell}\sigma_1(n)
-\frac{1}{24}\sigma_{\ell-1}\left(\frac{n}{N}\right)
\]
and
\[n\mapsto \sum_{k<n/N}\sigma_1(k)\sigma_{\ell-1}(n-kN)-\frac{B_\ell}{2\ell}\sigma_1\left(\frac{n}{N}\right)
-\frac{1}{24}\sigma_{\ell-1}(n)
\]
are linear combinations of the sets of functions
\[
\bigcup_{(\psi,\phi,t)\in A_{N,\ell+2}^*}\left\{
n\mapsto \sigma_{\ell+1}^{\psi,\phi}\left(\frac{n}{t}\right)
\right\},
\]
\[
\bigcup_{(\psi,\phi,t)\in A_{N,\ell}^*}\left\{
n\mapsto n\sigma_{\ell-1}^{\psi,\phi}\left(\frac{n}{t}\right)
\right\},
\]
\[
\bigcup_{\substack{(d,M)\in\N^* \\ dM\mid N}}
\bigcup_{f\in H_{\ell+2}^*[\Gamma_0(M)]}
\left\{
n\mapsto
\widehat{f}\left(\frac{n}{d}\right)
\right\},
\]
\[
\bigcup_{\substack{(d,M)\in\N^* \\ dM\mid N}}
\bigcup_{f\in H_{\ell}^*[\Gamma_0(M)]}
\left\{
n\mapsto
n\widehat{f}\left(\frac{n}{d}\right)
\right\}.
\]
\end{proposition}

The same allows to generalise theorem \ref{thm:deux}. If $b\geq 1$ is an integer, denote by $X(b)$ the
set of Dirichlet characters of modulus $b$. By orthogonality, we have
\[
S[a,b](n)=\frac{1}{\varphi(b)}\sum_{\chi\in X(b)}\overline{\chi(a)}
\sum_{m=1}^{n-1}\chi(m)\sigma_1(m)\sigma_1(n-m).
\]
It follows that the function to be considered is now
\[
\frac{1}{\varphi(b)}\sum_{\chi\in X(b)}\overline{\chi(a)}E_2\cdot E_2\otimes\chi.
\]
We restrict to $b$ squarefree so that the Gauss sum associates to any character of modulus $b$ is non vanishing.
For $N\geq 1$, let $\chi_N^{(0)}$ be the principal character of modulus $N$. 
For $\chi\in X(b)$, we define $A_{N,k,\chi}^*$ as $A_{N,k}^*$ except we replace condition \eqref{eq:extra} by
\[
\psi\phi=\chi_N^{(0)}\chi.
\]
Then, similarly to the proposition \ref{prop:son}, we have the following proposition.
\begin{proposition}
Let $b\geq 1$ squarefree and $a\in[0,b-1]$ be integers. Then the function
\[
n\mapsto S[a,b](n)
-\frac{1}{24}[\delta(b\mid a)+\delta(b\mid n-a)]\sigma_1(n)
\]
is a linear combination of the set of functions
\[
\bigcup_{\chi\in X(b)}
\bigcup_{(\psi,\phi,t)\in A_{N,4,\chi}^*}\left\{
n\mapsto \sigma_{3}^{\psi,\phi}\left(\frac{n}{t}\right)
\right\},
\]
\[
\bigcup_{\chi\in X(b)}
\bigcup_{(\psi,\phi,t)\in A_{N,2,\chi}^*}\left\{
n\mapsto n\sigma_{1}^{\psi,\phi}\left(\frac{n}{t}\right)
\right\},
\]
\[
\bigcup_{\chi\in X(b)}
\bigcup_{\substack{(d,M)\in\N^* \\ dM\mid N}}
\bigcup_{f\in H_{4}^*[\Gamma_0(M),\chi_N^{(0)}\chi]}
\left\{
n\mapsto
\widehat{f}\left(\frac{n}{d}\right)
\right\},
\]
\[
\bigcup_{\chi\in X(b)}
\bigcup_{\substack{(d,M)\in\N^* \\ dM\mid N}}
\bigcup_{f\in H_{2}^*[\Gamma_0(M),\chi_N^{(0)}\chi]}
\left\{
n\mapsto
n\widehat{f}\left(\frac{n}{d}\right)
\right\}
\]
\[
\{n\mapsto n\sigma_1(n)\}
\]
where $N$ is the least common multiple of $2$ and $b^2$ and $\delta(b\mid n-a)$ is $1$
if $n\equiv a\pmod{b}$ and $0$ otherwise.
\end{proposition}
\section{Convolutions of the divisor sum}
\subsection{Level 3}
By lemma \ref{lem:0040}, we have
\[
\qm{4}{2}[\Gamma_0(3)]
=
\m{4}[\Gamma_0(3)]\oplus D\m{2}[\Gamma_0(3)]\oplus \C DE_2.
\]
The vector space $\m{4}[\Gamma_0(3)]$ has dimension $2$
and is spanned by the two linearly independent forms $E_4$ and $E_{4,3}$.
The vector space $\m{2}[\Gamma_0(3)]$ has dimension $1$
and is spanned by $\Phi_{1,3}$.
Computing the first Fourier coefficients, we therefore find that
\begin{equation}\label{eq:H3}
H_3
=
\frac{1}{10}E_4
+\frac{9}{10}E_{4,3}
+4D\Phi_{1,3}
+4 DE_2.
\end{equation}

Comparing with the Fourier expansion in \eqref{eq:H3} leads to the corresponding result in theorem \ref{thm:un}.
\subsection{Level 5}
By lemma \ref{lem:0040}, we have
\[
\qm{4}{2}[\Gamma_0(5)]
=
\m{4}[\Gamma_0(5)]\oplus D\m{2}[\Gamma_0(5)]\oplus \C DE_2.
\]
The vector space $\m{4}[\Gamma_0(5)]$ has dimension $3$
and is spanned by the linearly independent forms $E_4$, $E_{4,5}$ and $\Delta_{4,5}$.
The vector space $\m{2}[\Gamma_0(5)]$ has dimension $1$
and is spanned by $\Phi_{1,5}$.
Computing the first Fourier coefficients, we therefore find that
\begin{equation}\label{eq:H5}
H_5
=
\frac{1}{26}E_4
+\frac{25}{26}E_{4,5}
-\frac{288}{65}\Delta_{4,5}
+\frac{24}{5}D\Phi_{1,5}
+\frac{12}{5}DE_2
\end{equation}

Comparing with the Fourier expansion in \eqref{eq:H5} leads to the corresponding result in theorem \ref{thm:un}.

\subsection{Level 6}
By lemma \ref{lem:0040}, we have
\[
\qm{4}{2}[\Gamma_0(6)]
=
\m{4}[\Gamma_0(6)]\oplus D\m{2}[\Gamma_0(6)]\oplus \C DE_2.
\]
The vector space $\m{4}[\Gamma_0(6)]$ has dimension $5$
and is spanned by the five linearly independent forms $E_4$, $E_{4,2}$ ,$E_{4,3}$ $E_{4,6}$ and $\Delta_{4,6}$.
The vector space $\m{2}[\Gamma_0(6)]$ has dimension $3$
and is spanned by the three linearly independent forms $\Phi_{1,2}$, $\Phi_{1,3}$ and $\Phi_{3,6}$.
Computing the first Fourier coefficients, we therefore find that
\begin{equation}\label{eq:H6}
H_6
=
\frac{1}{50}E_4
+\frac{2}{25}E_{4,2}
+\frac{9}{50}E_{4,3}
+\frac{18}{25}E_{4,6}
-\frac{24}{5}\Delta_{4,6}
+2D\Phi_{1,3}
+3D\Phi_{3,6}
+2DE_2.
\end{equation}

Comparing with the Fourier expansion in \eqref{eq:H6} leads to the corresponding result in theorem \ref{thm:un}.
\subsection{Level 7}

By lemma \ref{lem:0040}, we have
\[
\qm{4}{2}[\Gamma_0(7)]
=
\m{4}[\Gamma_0(7)]\oplus D\m{2}[\Gamma_0(7)]\oplus \C DE_2.
\]
The vector space $\m{4}[\Gamma_0(7)]$ has dimension $3$
and is spanned by the three linearly independent forms $E_4$, $E_{4,7}$ and $\Delta_{4,7}$.
The vector space $\m{2}[\Gamma_0(7)]$ has dimension $1$
and is spanned by the form $\Phi_{1,7}$.
Computing the first Fourier coefficients, we therefore find that
\begin{equation}\label{eq:H7}
H_7
=
\frac{1}{50}E_4
+\frac{49}{50}E_{4,7}
-\frac{288}{35}\Delta_{4,7}
+\frac{36}{7}D\Phi_{1,7}
+\frac{12}{7}DE_2.
\end{equation}
Comparing with the Fourier expansion in \eqref{eq:H7} leads to the corresponding result in theorem \ref{thm:un}.
\subsection{Level 8}

By lemma \ref{lem:0040}, we have
\[
\qm{4}{2}[\Gamma_0(8)]
=
\m{4}[\Gamma_0(8)]\oplus D\m{2}[\Gamma_0(8)]\oplus \C DE_2.
\]
The vector space $\m{4}[\Gamma_0(8)]$ has dimension $5$
and is spanned by the five linearly independent forms $E_4$, $E_{4,2}$, $E_{4,4}$, $E_{4,8}$ and $\Delta_{4,8}$.
The vector space $\m{2}[\Gamma_0(8)]$ has dimension $3$
and is spanned by the forms $\Phi_{1,4}$, $\Phi_{1,8}$ and 
\[\Phi_{1,4,2}\coloneqq z\mapsto\Phi_{1,4}(2z).\]
Computing the first Fourier coefficients, we therefore find that
\begin{equation}\label{eq:H8}
H_8
=
\frac{1}{80}E_4
+\frac{3}{80}E_{4,2}
+\frac{3}{20}E_{4,4}
+\frac{4}{5}E_{4,8}
-9\Delta_{4,8}
+\frac{21}{4}D\Phi_{1,8}
+\frac{3}{2}DE_2.
\end{equation}
Comparing with the Fourier expansion in \eqref{eq:H8} leads to the corresponding result in theorem \ref{thm:un}.
\subsection{Level 9}\label{cas9}
By lemma \ref{lem:0040}, we have
\[
\qm{4}{2}[\Gamma_0(9)]
=
\m{4}[\Gamma_0(9)]\oplus D\m{2}[\Gamma_0(9)]\oplus \C DE_2.
\]
The vector space $\m{4}[\Gamma_0(9)]$ has dimension $5$
and is spanned by the five linearly independent forms $E_4$, $E_{4}\otimes\chi_3$, $E_{4,3}$, $E_{4,9}$ and $\Delta_{4,9}$.
The vector space $\m{2}[\Gamma_0(9)]$ has dimension $3$
and is spanned by the forms $\Phi_{1,3}$, $\Phi_{1,3}\otimes\chi_3$ and $\Phi_{1,9}$.

Computing the first Fourier coefficients, we therefore find that
\begin{equation}\label{eq:H9}
H_9
=
\frac{1}{90}E_4
+\frac{4}{45}E_{4,3}
+\frac{9}{10}E_{4,9}
-\frac{32}{3}\Delta_{4,9}
+\frac{16}{3}D\Phi_{1,9}
+\frac{4}{3}DE_2.
\end{equation}
Comparing with the Fourier expansion in \eqref{eq:H9} leads to the corresponding result in theorem \ref{thm:un}.
\subsection{Level 10}\label{cas10}
By lemma \ref{lem:0040}, we have
\[
\qm{4}{2}[\Gamma_0(10)]
=
\m{4}[\Gamma_0(10)]\oplus D\m{2}[\Gamma_0(10)]\oplus \C DE_2.
\]
The vector space $\m{4}[\Gamma_0(10)]$ has dimension $7$
and is spanned by the seven linearly independent forms $E_4$, $E_{4,2}$, $E_{4,5}$, $E_{4,10}$, $\Delta_{4,10}$, $\Delta_{4,5}$
 and 
\[
F_{4,5,2}
\coloneqq
z\mapsto\Delta_{4,5}(2z).
\]
The vector space $\m{2}[\Gamma_0(10)]$ has dimension $3$
and is spanned by the forms $\Phi_{1,10}$, $\Phi_{1,5}$ and
\[
\Phi_{1,5,2}
\coloneqq
z\mapsto\Phi_{1,5}(2z).
\]

Computing the first Fourier coefficients, we therefore find that
\begin{multline}\label{eq:H10}
H_{10}
=
\frac{1}{130}E_4
+\frac{2}{65}E_{4,2}
+\frac{5}{26}E_{4,5}
+\frac{10}{13}E_{4,10}
-\frac{24}{5}\Delta_{4,10}
-\frac{432}{65}\Delta_{4,5}
-\frac{1728}{65}F_{4,5,2}\\
+\frac{27}{5}D\Phi_{1,10}
+\frac{6}{5}DE_2.
\end{multline}
Comparison with the Fourier expansion in \eqref{eq:H10} leads to the corresponding result in theorem \ref{thm:un}.
\subsection{Level 11}\label{cas11}
By lemma \ref{lem:0040}, we have
\[
\qm{4}{2}[\Gamma_0(11)]
=
\m{4}[\Gamma_0(11)]\oplus D\m{2}[\Gamma_0(11)]\oplus \C DE_2.
\]
The vector space $\m{4}[\Gamma_0(11)]$ has dimension $4$
and is spanned by the four linearly independent forms $E_4$, $E_{4,11}$, $\Delta_{4,11,1}$ and $\Delta_{4,11,2}$.
Let $F_1$ be the parabolic form of weight $4$ and level $11$ given by
\[
F_1(z)=[\Delta(z)\Delta(11z)]^{1/6}=e^{4\pi iz}-4e^{6\pi iz}+2e^{8\pi iz}+8e^{10\pi iz}+O(e^{12\pi iz}).
\]
Let $T_2$ be the Hecke operator of level $11$ given by
\[
T_2  : 
\sum_{m\in\Z}\widehat{f}(m)e^{2\pi imz}
\mapsto 
\sum_{m\in\Z}
\left[
\sum_{\substack{d\in\N\\ d\mid (m,2)\\ (d,11)=1}}d^{k-1}\widehat{f}\left(\frac{2m}{d^2}\right)
\right]
e^{2\pi imz}.
\]
It sends a parabolic form of weight $4$ and level $11$ to another one. Let
\[
F_2=T_2F_1=e^{2\pi iz}+2e^{4\pi iz}-5e^{6\pi iz}-2e^{8\pi iz}+9e^{10\pi iz}+O(e^{12\pi iz}).
\]
There exists $\lambda_1$ and $\lambda_2$ such that
\[
\Delta_{4,11,1}=F_2+\lambda_1F_1 \quad\text{ and }\quad
\Delta_{4,11,2}=F_2+\lambda_2F_1.
\]
For $j\in\{1,2\}$, it follows that
\[
\tau_{4,11,j}(2)=2+\lambda_j \text{ and } \tau_{4,11,j}(4)=-2+2\lambda_j.
\]
Since $\Delta_{4,11,j}$ is primitive, it satisfies \eqref{eq:hecke} hence $\lambda_j^2-2\lambda_j-2=0$. In other words
\[
X^2-2X-2=(X-\lambda_1)(X-\lambda_2).
\]
This provides a way to compute the Fourier coefficients of $\Delta_{4,11,1}$ and $\Delta_{4,11,2}$ from the ones of $\Delta$
and proves that these coefficients live in $\Q(t)$ where $t$ is a root of $X^2-2X-2$.

The vector space $\m{2}[\Gamma_0(11)]$ has dimension $2$
and is spanned by the form $\Phi_{1,11}$ and its unique primitive form 
\[
\Delta_{2,11}=[\Delta(z)\Delta(11z)]^{1/12}.
\]

\begin{table}[!ht]
\centering
\begin{tabular}{|c||c|c|c|c|c|c|c|c|c|c|c|}
\hline
$n$             & $1$  & $2$  & $3$   & $4$  & $5$  & $6$  & $7$  & $8$    & $9$    & $10$  & $11$
\\ \hline
$\tau_{2,11}(n)$ & $1$  & $-2$ & $-1$  & $2$ & $1$ & $2$  & $-2$ & $0$   & $-2$  & $-2$ & $1$
\\ \hline
\hline
$n$             & $12$ & $13$ & $14$  & $15$ & $16$ & $17$ & $18$ & $19$   & $20$   & $21$  & $22$ 
\\ \hline
$\tau_{2,11}(n)$ & $-2$ & $4$ & $4$ & $-1$  & $-4$ & $-2$ & $4$ & $0$ & $2$ & $2$  & $-2$ 
\\ \hline
\end{tabular}
\caption{First Fourier coefficients of $\Delta_{2,11}$}
\end{table}

Computing the first Fourier coefficients, we therefore find that
\begin{multline}\label{eq:H11}
H_{11}
=
\frac{1}{122}E_4
+\frac{121}{122}E_{4,11}
-\frac{192t+4128}{671}\Delta_{4,11,1}
+\frac{192t-4512}{671}\Delta_{4,11,2}
\\
+\frac{60}{11}D\Phi_{1,11}
+\frac{12}{11}DE_2.
\end{multline}
Comparison with the Fourier expansion in \eqref{eq:H11} leads to theorem \ref{thm:onze}.
\subsection{Level 13}\label{cas13}
By lemma \ref{lem:0040}, we have
\[
\qm{4}{2}[\Gamma_0(13)]
=
\m{4}[\Gamma_0(13)]\oplus D\m{2}[\Gamma_0(13)]\oplus \C DE_2.
\]
The vector space $\m{4}[\Gamma_0(13)]$ has dimension $5$
and is spanned by the five linearly independent forms $E_4$, $E_{4,13}$, $\Delta_{4,13,1}$, $\Delta_{4,13,2}$ and
$\Delta_{4,13,3}.$
The vector space $\m{2}[\Gamma_0(13)]$ has dimension $1$
and is spanned by the form $\Phi_{1,13}$.

\begin{table}[!ht]
\centering
\begin{tabular}{|c||c|c|c|c|c|c|c|c|c|c|c|}
\hline
$n$             & $1$  & $2$  & $3$   & $4$  & $5$  & $6$  & $7$  & $8$    & $9$    & $10$  & $11$
\\ \hline
$\tau_{4,13,1}(n)$ & $1$  & $-5$ & $-7$  & $17$ & $-7$ & $35$  & $-13$ & $-45$   & $22$  & $35$ & $-26$
\\ \hline
\hline
$n$             & $12$ & $13$ & $14$  & $15$ & $16$ & $17$ & $18$ & $19$   & $20$   & $21$  & $22$ 
\\ \hline
$\tau_{4,13,1}(n)$ & $-119$ & $13$ & $65$ & $49$  & $89$ & $77$ & $-110$ & $-126$ & $-119$ & $91$  & $130$ 
\\ \hline
\end{tabular}
\caption{First Fourier coefficients of $\Delta_{4,13,1}$}
\end{table}

Computing the first Fourier coefficients, we therefore find that
\begin{multline}\label{eq:H13}
H_{13}
=
\frac{1}{170}E_4
+\frac{169}{170}E_{4,13}
+\frac{288u-1728}{221}\Delta_{4,13,2}
-\frac{288u+1440}{221}\Delta_{4,13,3}
\\
+\frac{72}{13}D\Phi_{1,13}
+\frac{12}{13}DE_2.
\end{multline}
Comparison with the Fourier expansion in \eqref{eq:H13} leads to the theorem \ref{thm:treize}.
\subsection{Level 14}\label{cas14}
By lemma \ref{lem:0040}, we have
\[
\qm{4}{2}[\Gamma_0(14)]
=
\m{4}[\Gamma_0(14)]\oplus D\m{2}[\Gamma_0(14)]\oplus \C DE_2.
\]
The vector space $\m{4}[\Gamma_0(14)]$ has dimension $8$
and is spanned by the eight linearly independent forms $E_4$, $E_{4,2}$, $E_{4,7}$, $E_{4,14}$, $\Delta_{4,7}$,
\[
F_{4,7,2} \colon z \mapsto \Delta_{4,7}(2z),
\]
$\Delta_{4,14,1}$ and $\Delta_{4,14,2}$. Another basis of the subspace of parabolic forms
is $\Delta_{4,7}$, $F_{4,7,2}$, $\Delta_{2,14}^2$ and $\Delta_{2,14}\Phi_{1,14}$ where 
\[
\Delta_{2,14}(z)\coloneqq[\Delta(z)\Delta(2z)\Delta(7z)\Delta(14z)]^{1/24}
\] 
is the unique primitive form of weight $2$ on $\Gamma_0(14)$. We echelonise this second basis by defining
\begin{align*}
J_1 &= -\frac{11}{28}\Delta_{4,7}-\frac{22}{7}F_{4,7,2}+\frac{11}{7}\Delta_{2,14}^2+\frac{39}{28}\Delta_{2,14}\Phi_{1,14}
&= e^{2\pi iz}+O\left(e^{10\pi iz}\right)\\
J_2 &= -\frac{13}{56}\Delta_{4,7}+\frac{1}{7}F_{4,7,2}+\frac{3}{7}\Delta_{2,14}^2+\frac{13}{56}\Delta_{2,14}\Phi_{1,14}
&= e^{4\pi iz}+O\left(e^{10\pi iz}\right)\\
J_3 &= \frac{13}{56}\Delta_{4,7}+\frac{19}{14}F_{4,7,2}-\frac{13}{14}\Delta_{2,14}^2-\frac{13}{56}\Delta_{2,14}\Phi_{1,14}
&= e^{6\pi iz}+O\left(e^{10\pi iz}\right)\\
J_4 &= -\frac{13}{56}\Delta_{4,7}-\frac{6}{7}F_{4,7,2}+\frac{3}{7}\Delta_{2,14}^2+\frac{13}{56}\Delta_{2,14}\Phi_{1,14}
&= e^{8\pi iz}+O\left(e^{10\pi iz}\right).
\end{align*}
We then have
\[
\Delta_{4,14,j}=J_1+b_jJ_2+c_jJ_3+d_jJ_4.
\]
From $\tau_{4,14,j}(4)=\tau_{4,14,j}(2)^2$ we deduce $d_j=b_j^2$. Then, from $\tau_{4,14,j}(6)=\tau_{4,14,j}(2)\tau_{4,14,j}(3)$ and
$\tau_{4,14,j}(8)=\tau_{4,14,j}(2)\tau_{4,14,j}(4)$ we respectively deduce
\begin{align*}
2b_j+b_jc_j+2c_j &= -4\\
b_j^3-b_j^2+6b_j+4c_j &= 8
\end{align*}
that is  
\[
c_j=-\frac{1}{4}b_j^3+\frac{1}{4}b_j^2-\frac{3}{2}b_j+2
\]
and
\[
(b_j-2)(b_j+2)(b_j^2+b_j+8)=0.
\]
Since the coefficients of $\Delta_{4,14,j}$ are all totally real, we must have
\begin{align*}
\Delta_{4,14,1} &= J_1+2J_2-2J_3+4J_4\\
\Delta_{4,14,2} &= J_1-2J_2+8J_3+4J_4.
\end{align*}
Finally,
\begin{align}
\label{eq:4141}
\Delta_{4,14,1} &= -\frac{9}{4}\Delta_{4,7}-9F_{4,7,2}+6\Delta_{2,14}^2+\frac{13}{4}\Delta_{2,14}\Phi_{1,14}\\
\label{4142}
\Delta_{4,14,2} &=\Delta_{4,7}+4F_{4,7,2}-5\Delta_{2,14}^2.
\end{align}
Equations \eqref{eq:4141} and \eqref{4142} allow to compute the first terms of the sequences $\tau_{4,14,1}$ and $\tau_{4,14,2}$.
The vector space $\m{2}[\Gamma_0(14)]$ has dimension $4$
and is spanned by the forms $\Phi_{1,7}$, $\Phi_{1,14}$, $\Phi_{2,14}$ and its unique primitive form $\Delta_{2,14}$.

\begin{table}[!ht]
\centering
\begin{tabular}{|c||c|c|c|c|c|c|c|c|c|c|c|}
\hline
$n$             & $1$  & $2$  & $3$   & $4$  & $5$  & $6$  & $7$  & $8$    & $9$    & $10$  & $11$
\\ \hline
$\tau_{2,14}(n)$  & $1$ & $-1$ & $-2$ & $1$ & $0$ & $2$ & $1$ & $-1$ & $1$ & $0$ & $0$
\\ \hline
\hline
$n$             & $12$ & $13$ & $14$  & $15$ & $16$ & $17$ & $18$ & $19$   & $20$   & $21$  & $22$ 
\\ \hline
$\tau_{2,14}(n)$ & $-2$ & $-4$ & $-1$ & $0$ & $1$ & $6$ & $-1$ & $2$ & $0$ & $-2$ & $0$
\\ \hline
\end{tabular}
\caption{First Fourier coefficients of $\Delta_{2,14}$}
\end{table}

Computing the first Fourier coefficients, we therefore find that
\begin{multline}\label{eq:H14}
H_{14}
=
\frac{1}{250}E_4
+\frac{2}{125}E_{4,2}
+\frac{49}{250}E_{4,7}
+\frac{98}{125}E_{4,14}
-\frac{864}{175}\Delta_{4,7}
-\frac{3456}{175}F_{4,7,2}\\
-\frac{48}{7}\Delta_{4,14,1}
-\frac{72}{25}\Delta_{4,14,2}
+\frac{39}{7}D\Phi_{1,14}
+\frac{6}{7}DE_2.
\end{multline}
Comparison with the Fourier expansion in \eqref{eq:H14} leads to theorem \ref{thm:quatorze}.
\subsection{Convolutions of level 1, 2, 4}

The convolutions of level dividing $4$ were evaluated in \cite[Proposition 7]{GT0509205}. We obtained
\[
W_1(n)
=
\frac{5}{12}\sigma_3(n)
-\frac{n}{2}\sigma_1(n)
+\frac{1}{12}\sigma_1(n)
\]
from the equality
\begin{equation*}
E_2^2=E_4+12DE_2
\end{equation*}
in $\qm{4}{2}[\Gamma_0(1)]$;
\[
W_2(n)
=
\frac{1}{12}\sigma_3(n)
+\frac{1}{3}\sigma_3\left(\frac{n}{2}\right)
-\frac{1}{8}n\sigma_1(n)
-\frac{1}{4}n\sigma_1\left(\frac{n}{2}\right)
+\frac{1}{24}\sigma_1(n)
+\frac{1}{24}\sigma_1\left(\frac{n}{2}\right)
\]
from the equality
\[
H_2=\frac{1}{5}E_4+\frac{4}{5}E_{4,2}+3D\Phi_{1,2}+6DE_2
\]
in $\qm{4}{2}[\Gamma_0(2)]$; and
\begin{multline*}
W_4(n)
=
\frac{1}{48}\sigma_3(n)
+\frac{1}{16}\sigma_3\left(\frac{n}{2}\right)
+\frac{1}{3}\sigma_3\left(\frac{n}{4}\right)
-\frac{1}{16}n\sigma_1(n)
-\frac{1}{4}n\sigma_1\left(\frac{n}{4}\right)\\
+\frac{1}{24}\sigma_1(n)
+\frac{1}{24}\sigma_1\left(\frac{n}{4}\right).
\end{multline*}
from the equality
\[
H_4
=
\frac{1}{20}E_4
+
\frac{3}{20}E_{4,2}
+\frac{4}{5}E_{4,4}
+\frac{9}{2}D\Phi_{1,4}
+3DE_2.
\]
in $\qm{4}{2}[\Gamma_0(4)]$.
\section{Twisted convolution sums}

Let $\chi^{(0)}_3$ be the principal character of modulus $3$. Remarking
that
\[
S[0,3](n)=\sum_{a+b=n}\sigma_1(a)\sigma_1(b)-\sum_{a+b=n}\chi^{(0)}_3(a)\sigma_1(a)\sigma_1(b),
\]
we consider $E_2^2-E_2(E_2\otimes\chi^{(0)}_3)$. Since 
$E_2\otimes\chi^{(0)}_3\in\qm{2}{1}[\Gamma_0(9),\chi^{(0)}_3]=\qm{2}{1}[\Gamma_0(9)]$,
we have
\[
E_2^2-E_2(E_2\otimes\chi^{(0)}_3)\in\qm{4}{2}[\Gamma_0(9)].
\]
We use the same method and notations as in \S\ref{cas9}. We compute
\[
E_2^2-E_2(E_2\otimes\chi^{(0)}_3)
=
\frac{11}{30}E_4+\frac{10}{3}E_{4,3}-\frac{27}{10}E_{4,9}+32\Delta_{4,9}
+16D\Phi_{1,3}-16D\Phi_{1,9}
+12DE_2.
\]
The evaluation of $S[1,3]$ given in theorem \ref{thm:deux} follows by comparison of the Fourier expansions.

We compute $S[1,3]$ after having remarked that
\[
\frac{\chi_3^{(0)}(n)+\chi_3(n)}{2}
=
\begin{cases}
1 &\text{if $n\equiv 1\pmod{3}$}\\
0 &\text{otherwise.}
\end{cases}
\]
Hence, the function to be linearised here is
\[
\frac{1}{2}E_2[E_2\otimes\chi_3^{(0)}+E_2\otimes\chi_3]
\]
whose $n$th Fourier coefficient ($n\in\N^*$) is
\[
-24\delta(3\mid n-1)\sigma_1(n)+576S[1,3](n).
\]
This is again a quasimodular form in $\qm{4}{2}[\Gamma_0(9)]$, and as in \S\ref{cas9}, we linearise it as
\[
\frac{19}{60}E_4
+\frac{1}{20}E_4\otimes\chi_3
-\frac{5}{3}E_{4,3}
+\frac{27}{20}E_{4,9}
+32\Delta_{4,9}
-8D\Phi_{1,3}
-6D(\Phi_{1,3}\otimes\chi_3)
+8D\Phi_{1,9}.
\]
The evaluation of $S[1,3]$ given in theorem \ref{thm:deux} follows by comparison of the Fourier expansions.

The evaluation of $S[2,3]$ follows immediately from
\[
S[0,3](n)+S[1,3](n)+S[2,3](n)=W_1(n)
\]
and theorem \ref{thm:un}.
\section{On identities by Melfi}
The first three identities of theorem \ref{thm:trois} are a direct consequence of the following ones:
\[
E_2E_4\in\qm{6}{1}[\Gamma_0(1)]
=
\m{6}[\Gamma_0(1)]\oplus D\m{4}[\Gamma_0(1)]
=
\C E_6\oplus\C DE_4,
\]
and
\begin{align*}
E_2E_{4,2},\, E_4E_{2,2}\in\qm{6}{1}[\Gamma_0(2)]
&=
\m{6}[\Gamma_0(2)]\oplus D\m{4}[\Gamma_0(2)]
\\
&=
\C E_6\oplus\C E_{6,2}\oplus\C DE_4\oplus\C DE_{4,2}
\end{align*}
which imply by comparison of the first Fourier coefficients
\[
E_2E_4=E_6+3DE_4,
\]
\[
E_2E_{4,2}
=
\frac{1}{21}E_6+\frac{20}{21}E_{6,2}+3D E_{4,2}
\]
and
\[
E_4E_{2,2}
=
\frac{5}{21}E_6+\frac{16}{21}E_{6,2}+\frac{3}{2}D E_{4}.
\]

The last three identities of theorem \ref{thm:trois} are a direct consequence of the following ones:
\[
E_2E_6\in\qm{8}{1}[\Gamma_0(1)]
=
\m{8}[\Gamma_0(1)]\oplus D\m{6}[\Gamma_0(1)]
=
\C E_8\oplus\C DE_6,
\]
and
\begin{align*}
E_{2,2}E_{6},\, E_2E_{6,2}\in\qm{8}{1}[\Gamma_0(2)]
&=
\m{8}[\Gamma_0(2)]\oplus D\m{6}[\Gamma_0(2)]
\\
&=
\C E_8\oplus\C E_{8,2}\oplus\C\Delta_{8,2}\oplus\C DE_6\oplus\C DE_{6,2}
\end{align*}
which imply by comparison of the first Fourier coefficients
\[
E_2E_6=E_8+2DE_6,
\]
\[
E_{2,2}E_6
=
\frac{21}{85}E_8+\frac{64}{85}E_{8,2}-\frac{2016}{17}\Delta_{8,2}+DE_6
\]
and
\[
E_2E_{6,2}
=
\frac{1}{85}E_8+\frac{84}{85}E_{8,2}-\frac{504}{17}\Delta_{8,2}+2DE_{6,2}.
\]
\section{On some identities of Lahiri type}
\subsection{Method}

For $\vec{a}\coloneqq (a_1,\dotsc,a_r)\in\N^r$, $\vec{b}\coloneqq(b_1,\dotsc,b_r)\in(2\N+1)^r$ and
$\vec{N}\coloneqq (N_1,\dotsc,N_r)\in\N^{*r}$
the sum $S[\vec{a},\vec{b},\vec{N}]$ defined in \eqref{eq:defS} is relied to the
quasimodular forms \textit{via} the function
\begin{equation}\label{eq:sela}
D^{a_1}E_{b_1+1,N_1}\dotsm D^{a_r}E_{b_r+1,N_r}
\in
\qm{b_1+\dotsm+b_r+r+2(a_1+\dotsc+a_r)}{a_1+\dotsc+a_r+t(\vec{b})}\left[\Gamma_0\left(\lcm(N_1,\dotsc,N_r)\right)\right]
\end{equation}
where
\[
t(\vec{b})
=
\#\{i\in\{1,\dotsc,r\} \colon b_i=1\}.
\]
Since we always can consider that the coordinates of $\vec{a}$ are given in increasing order, let
$\ell$ be the nonnegative integer such that $a_1=\dotsm=a_\ell=0$ and $a_{\ell+1}\neq 0$ (we take $\ell=0$
if $\vec{a}$ has all its coordinates positive). We consider the function
\begin{align*}
\Psi_{\vec{a},\vec{b},\vec{N}}&\coloneqq \prod_{j=1}^\ell(E_{b_j+1,N_1}-1)\prod_{j=\ell+1}^rD^{a_j}E_{b_j+1,N_j}\\
& \in
\bigoplus_{k=b_{j+1}+\dotsm+b_r+r+2(a_{j+1}+\dotsc+a_r)-j}^{b_1+\dotsm+b_r+r+2(a_{j+1}+\dotsc+a_r)}
\qm{k}{a_1+\dotsc+a_r+t(\vec{b})}\left[\Gamma_0\left(\lcm(N_1,\dotsc,N_r)\right)\right]
\end{align*} 
We have
\[
\Psi_{\vec{a},\vec{b},\vec{N}}(z)
=
(-2)^r\left[\prod_{j=1}^{r}\frac{b_j+1}{B_{b_j+1}}\right]\sum_{n=1}^{+\infty}S[\vec{a},\vec{b},\vec{N}]e^{2\pi inz}.
\]

The evaluation of $S[(0,1,1),(1,1,1)]$ is a consequence (by lemma \ref{lem:0040}) of
\begin{multline*}
(E_1-1)(DE_2)^2
\in
\\
\C E_8\oplus\C DE_6\oplus\C D^2E_4\oplus\C D^3E_2\oplus\C E_{10}\oplus\C DE_8\oplus\C D^2E_6\oplus\C D^3E_4\oplus\C D^4E_2.
\end{multline*}
The comparison of the first Fourier coefficients leads to
\[
(E_2-1)(DE_2)^2
=
-\frac{1}{5}D^2E_4-2D^3E_2+\frac{2}{21}D^2E_6+\frac{4}{5}D^3E_4+6D^4E_2.
\]
Hence the evaluation of $S[(0,1,1),(1,1,1)]$ given in theorem \ref{thm:quatre}.

The evaluation of $S[(0,0,0,1,1),(1,1,1,1,1)]$ is a consequence (by lemma \ref{lem:0040}) of
\[
(E_2-1)^3(DE_2)^2
\in
\C\Delta\oplus\C D\Delta\bigoplus_{i=1}^7\bigoplus_{j=0}^{7-i}\C D^j E_{2i}.
\]
The comparison of the first Fourier coefficients leads to
\begin{multline*}
(E_2-1)^3(DE_2)^2
=
-\frac{8}{35}\Delta+\frac{8}{35}D\Delta
-2D^3E_2+18D^4E_2-\frac{216}{5}D^5E_2+\frac{144}{5}D^6E_2
\\
-\frac{1}{5}D^2E_4+\frac{12}{5}D^3E_4-\frac{234}{35}D^4E_4+\frac{171}{35}D^5E_4
\\
+\frac{2}{7}D^2E_6-\frac{9}{7}D^3E_6+\frac{25}{21}D^4E_6
-\frac{1}{6}D^2E_8+\frac{4}{15}D^3E_8
+\frac{2}{55}D^2E_{10}.
\end{multline*}
Hence the evaluation of $S[(0,0,0,1,1),(1,1,1,1,1)]$ given in theorem \ref{thm:quatre}.

We leave the proofs of propositions \ref{prop:1721} and \ref{prop:1715} to the reader. They are obtained 
from the linearisations of
\begin{equation}\label{eq:point}
(E_{2,2}-1)DE_{2,5}\in\qm{4}{2}[\Gamma_0(5)]\oplus \qm{6}{3}[\Gamma_0(10)]
\end{equation}
and
\[
DE_2DE_{2,5}\in\qm{8}{5}[\Gamma_0(5)].
\]

\subsection{Primitive forms of weight 6 and level 5 or 10}\label{pa:prec} 

For the evaluation of \eqref{eq:point} we remark that $\Delta_{4,5}\Phi_{1,5}$
is a parabolic modular form of weight $6$ and level $5$. Since the dimension of these forms is $1$, we have
\begin{equation}\label{eq:dsc}
\Delta_{6,5}(z)
=
\frac{1}{4}\left[\Delta(z)\Delta(5z)\right]^{1/6}\left[5E_2(5Z)-E_2(Z)\right].
\end{equation}
Equation \eqref{eq:dsc} provides a way to compute the few needed values of
$\tau_{6,5}$.

We also give expressions for $\Delta_{6,10,i}$ where $i\in\{1,2,3\}$. 
We shall use the second Hecke operator of level $10$ given by
\[
T_2  : 
\sum_{m\in\Z}\widehat{f}(m)e^{2\pi imz}
\mapsto 
\sum_{m\in\Z}
\left[
\sum_{\substack{d\in\N\\ d\mid (m,2)\\ (d,10)=1}}d^{k-1}\widehat{f}\left(\frac{2m}{d^2}\right)
\right]
e^{2\pi imz}
=
\sum_{m\in\Z}
\widehat{f}(2m)
e^{2\pi imz}.
\]
The space of parabolic modular forms of weight $6$ and level $10$ has dimension $5$. A basis is given by
\begin{align*}
\Delta_{6,5}(z) &= \left[\Delta(z)\Delta(5z)\right]^{1/6}\Phi_{1,5}(z)\\
F_{6,5,2}(z) &= \Delta_{6,5}(2z)\\
F(z) &= 3\left[\Delta(z)\Delta(5z)\right]^{1/6}\Phi_{1,10}(z)\\
F_1(z) &= \left[\Delta(z)\Delta(5z)\right]^{1/6}\Phi_{1,2}(z)\\
F_2(z) &= T_2F(z).
\end{align*}
To simplify the computations, we use an echelonised basis:
\begin{align*}
V_1&=
-\frac{4}{15}\Delta_{6,5}+\frac{31}{10}F_{6,5,2}+\frac{15}{32}F+\frac{1}{96}F_1+\frac{3}{80}F_2\\
&= e(z)+O(e(6z))\\
V_2&=
\frac{1}{20}\Delta_{6,5}+\frac{6}{5}F_{6,5,2}+\frac{1}{80}F_2\\
&= e(2z)+O(e(6z))\\
V_3&=
-\frac{1}{30}\Delta_{6,5}+\frac{7}{10}F_{6,5,2}+\frac{1}{32}F-\frac{1}{96}F_1+\frac{1}{80}F_2\\
&= e(3z)+O(e(6z))\\
V_4&=
-\frac{1}{40}\Delta_{6,5}-\frac{1}{10}F_{6,5,2}-\frac{1}{160}F_2\\
&= e(4z)+O(e(6z))\\
V_5&=
\frac{1}{75}\Delta_{6,5}-\frac{11}{50}F_{6,5,2}-\frac{11}{800}F-\frac{1}{480}F_1-\frac{3}{400}F_2\\
&= e(5z)+O(e(6z)).
\end{align*}
We deduce
\[
\Delta_{6,10,i}=V_1+b_iV_2+c_iV_3+d_iV_4+e_iV_5
\]
since $\tau_{6,10,i}(1)=1$. Now, from  $\tau_{6,10,i}(4)=\tau_{6,10,i}(2)^2$, we get
$d_i=b_i^2$ so that
\[
\Delta_{6,10,i}=V_1+b_iV_2+c_iV_3+b_i^2V_4+e_iV_5.
\]
Next, from $\tau_{6,10,i}(6)=\tau_{6,10,i}(2)\tau_{6,10,i}(3)$, and
$\tau_{6,10,i}(8)=\tau_{6,10,i}(2)\tau_{6,10,i}(4)$, we respectively get
\begin{equation}\label{eq:leci}
2c_i+b_ic_i+2e_i
=
10-2b_i-b_i^2
\end{equation}
and 
\begin{equation}\label{eq:nema}
8c_i-8e_i
=
24+16b_i+6b_i^2+b_i^3.
\end{equation}
Equations \eqref{eq:leci} and \eqref{eq:nema} give either $b_i=-4$ or $b_i\neq-4$ and
\[
c_i=4-\frac{1}{2}b_i+\frac{1}{4}b_i^2
\]
and
\[
e_i
=
1-\frac{5}{2}b_i-\frac{1}{2}b_i^2-\frac{1}{8}b_i^3.
\]
We first deal with the case $b_i\neq -4$. We then get
\[
\Delta_{6,10,i}
=
V_1+b_iV_2
+
\left(4-\frac{1}{2}b_i+\frac{1}{4}b_i^2\right)V_3
+b_i^2V_4
+\left(1-\frac{5}{2}b_i-\frac{1}{2}b_i^2-\frac{1}{8}b_i^3\right)V_5.
\]
Using $\tau_{6,10,i}(10)=\tau_{6,10,i}(2)\tau_{6,10,i}(5)$, we obtain
\begin{multline*}
b_i\left(1-\frac{5}{2}b_i-\frac{1}{2}b_i^2-\frac{1}{8}b_i^3\right)
=
\\
-30+15b_i-10\left(4-\frac{1}{2}b_i+\frac{1}{4}b_i^2\right)
+5b_i^2
+6\left(1-\frac{5}{2}b_i-\frac{1}{2}b_i^2-\frac{1}{8}b_i^3\right)
\end{multline*}
from what we get
\[
b_i\in\left\{-4,4,1-i\sqrt{31},1+i\sqrt{31}\right\}.
\]
The solution $b_i=-4$ is in this case not allowed whereas the solution $1\pm i\sqrt{31}$ are not possible since the coefficients
of a primitive form are totally real algebraic numbers. We thus obtain a first primitive form:
\[
\Delta_{6,10,1}
=
V_1+4V_2+6V_3+16V_4-25V_5.
\]
We assume now that $b_i=-4$ so that
\[
\Delta_{6,10,i}=V_1-4V_2+c_iV_3+16V_4+(c_i+1)V_5.
\]
From $\tau_{6,10,i}(15)=\tau_{6,10,i}(3)\tau_{6,10,i}(5)$, we obtain $c_i=24$ or $c_i=-26$.
We hence get the two other primitive forms
\[
\Delta_{6,10,2}
=
V_1-4V_2+24V_3+16V_4+25V_5
\]
and
\[
\Delta_{6,10,3}
=
V_1-4V_2-26V_3+16V_4-25V_5.
\]
We deduce the following expressions:
\begin{align}
\label{eq:dss1}
\Delta_{6,10,1}&=-\Delta_{6,5}+16F_{6,5,2}+F+\frac{1}{4}F_2\\
\label{eq:dss2}
\Delta_{6,10,2}&=-\frac{4}{3}\Delta_{6,5}+8F_{6,5,2}+\frac{7}{8}F-\frac{7}{24}F_1\\
\label{eq:dss3}
\Delta_{6,10,3}&=-\frac{1}{3}\Delta_{6,5}-16F_{6,5,2}+\frac{1}{3}F_1-\frac{1}{4}F_2.
\end{align}
Equations \eqref{eq:dss1} to \eqref{eq:dss3} provide a way to compute the few needed values of
$\tau_{6,10,i}$ for $i\in\{1,2,3\}$.

\subsection{Primitive forms of weight 8 and level 5}\label{pa:eight}

The method is the same as in \S\ref{pa:prec} so we will be more brief. The
space of parabolic forms of weight $8$ and level $5$ has dimension $3$ and a
basis is
\begin{align*}
G_1(z) &= [\Delta(z)\Delta(5z)]^{1/3}\\
G_2(z) &= [\Delta(z)\Delta(5z)]^{1/6}\Phi_{1,5}(z)^2\\
G_3 &= -\frac{1}{24}\left[E_4,\Phi_{1,2}\right]_1
\end{align*}
where $\left[\phantom{f},\phantom{g}\right]_1$ is the Rankin-Cohen bracket
here defined by
\[
\left[E_4,\Phi_{1,2}\right]_1
=
\frac{1}{2\pi i }(4E_4\Phi_{1,5}'-2E_4'\Phi_{1,5})
\]
(see \cite[part 1, \S E]{MR1221103} or \cite[partie I, \S 6]{MR2186573} for
more details). We echelonise this basis by defining:
\begin{align*}
W_1 &= \frac{46}{25}G_1+\frac{82}{25}G_2-\frac{3}{25}G_3
    &= e(z)+O(e(4z))\\
W_2 &= \frac{47}{375}G_1-\frac{76}{375}G_2+\frac{4}{375}G_3
    &= e(2z)+O(e(4z))\\
W_3 &= -\frac{41}{375}G_1-\frac{19}{750}G_2+\frac{1}{750}G_3
    &= e(3z)+O(e(4z)).
\end{align*}
The primitive forms are then
\[
\Delta_{8,5,i}=W_1+b_iW_2+c_iW_3.
\]
From $\tau_{8,5,i}(4)=\tau_{8,5,i}(2)^2-2^7$ and
 $\tau_{8,5,i}(6)=\tau_{8,5,i}(2)\tau_{8,5,i}(3)$
we get
\begin{align*}
c_i &= 78+2b_i-\frac{1}{2}b_i^2\\
(b_i+14)(b_i^2-20b_i+24)  &=0.
\end{align*}
Finally, defining $v$ as one of the roots of $X^2-20X+24$, we get
\begin{align}
\label{eq:est1}
\Delta_{8,5,1} &= \frac{16}{3}G_1+\frac{22}{3}G_2-\frac{1}{3}G_3\\
\label{eq:est2}
\Delta_{8,5,2} &= (12-v)G_1+G_2\\
\label{eq:est3}
\Delta_{8,5,3} &= (v-8)G_1+G_2.
\end{align}
Equations \eqref{eq:est1} to \eqref{eq:est3} provide a way to compute the few needed values of
$\tau_{8,5,i}$ for $i\in\{1,2,3\}$.

\bibliographystyle{alpha}
\bibliography{Royer_IJNT_06_biblio}

\begin{thebibliography}{HOSW02}

\bibitem[AAWa]{AAW16}
A.~Alaca, S.~Alaca, and K.S. Williams.
\newblock The convolution sum $\sum_{m<n/16}\sigma(m)\sigma(n-16m)$.
\newblock {\em Canad. Math. Bull.}
\newblock To appear.

\bibitem[AAWb]{AAW24}
A.~Alaca, S.~Alaca, and K.S. Williams.
\newblock The convolution sums $\sum_{\ell+24m=n}\sigma(\ell)\sigma(m)$ and
  $\sum_{3\ell+8m=n}\sigma(\ell)\sigma(m)$.
\newblock Preprint.

\bibitem[AAWc]{AAW12}
A.~Alaca, S.~Alaca, and K.S. Williams.
\newblock Evaluation of the convolution sums
  $\sum_{\ell+12m=n}\sigma(\ell)\sigma(m)$ and
  $\sum_{3\ell+4m=n}\sigma(\ell)\sigma(m)$.
\newblock {\em Advances in Theoretical and Applied Mathematics}.
\newblock To appear.

\bibitem[AAWd]{AAW18}
A.~Alaca, S.~Alaca, and K.S. Williams.
\newblock Evaluation of the convolution sums
  $\sum_{\ell+18m=n}\sigma(\ell)\sigma(m)$ and
  $\sum_{2\ell+9m=n}\sigma(\ell)\sigma(m)$.
\newblock Preprint.

\bibitem[AW]{AlWi}
S.~Alaca and K.S. Williams.
\newblock Evaluation of the convolution sums
  $\sum_{\ell+6m=n}\sigma(\ell)\sigma(m)$ and
  $\sum_{2\ell+3m=n}\sigma(\ell)\sigma(m)$.
\newblock Preprint.

\bibitem[BCP97]{MR1484478}
Wieb Bosma, John Cannon, and Catherine Playoust.
\newblock The {M}agma algebra system. {I}. {T}he user language.
\newblock {\em J. Symbolic Comput.}, 24(3-4):235--265, 1997.
\newblock Computational algebra and number theory (London, 1993).

\bibitem[Bes62]{Bes62}
M.~Besge.
\newblock {Extrait d'une lettre de M. Besge \`a M. Liouville}.
\newblock {\em J. Math. Pures Appl}, 7:256, 1862.

\bibitem[CW05]{CW05}
N.~Cheng and K.S. Williams.
\newblock Evaluation of some convolution sums involving the sum of divisors
  function.
\newblock {\em Yokohama Math. J.}, 52:39--57, 2005.

\bibitem[DI95]{MR1357209}
Fred Diamond and John Im.
\newblock Modular forms and modular curves.
\newblock In {\em Seminar on Fermat's Last Theorem (Toronto, ON, 1993--1994)},
  volume~17 of {\em CMS Conf. Proc.}, pages 39--133. Amer. Math. Soc.,
  Providence, RI, 1995.

\bibitem[DS05]{MR2112196}
Fred Diamond and Jerry Shurman.
\newblock {\em A first course in modular forms}, volume 228 of {\em Graduate
  Texts in Mathematics}.
\newblock Springer-Verlag, New York, 2005.

\bibitem[Gla85]{Gla84}
J.~W.~L. Glaisher.
\newblock On the square of the series in which the coefficients are the sums of
  the divisors of the exponents.
\newblock {\em Mess. Math.}, XV:156--163, 1885.

\bibitem[HOSW02]{MR1956253}
James~G. Huard, Zhiming~M. Ou, Blair~K. Spearman, and Kenneth~S. Williams.
\newblock Elementary evaluation of certain convolution sums involving divisor
  functions.
\newblock In {\em Number theory for the millennium, II (Urbana, IL, 2000)},
  pages 229--274. A K Peters, Natick, MA, 2002.

\bibitem[Koi84]{MR759465}
Masao Koike.
\newblock On {M}c{K}ay's conjecture.
\newblock {\em Nagoya Math. J.}, 95:85--89, 1984.

\bibitem[KZ95]{MR1363056}
Masanobu Kaneko and Don Zagier.
\newblock A generalized {J}acobi theta function and quasimodular forms.
\newblock In {\em The moduli space of curves (Texel Island, 1994)}, volume 129
  of {\em Progr. Math.}, pages 165--172. Birkh\"auser Boston, Boston, MA, 1995.

\bibitem[Lah46]{MR0020591}
D.~B. Lahiri.
\newblock On {R}amanujan's function {$\tau(n)$} and the divisor function
  {$\sigma\sb k(n)$}. {I}.
\newblock {\em Bull. Calcutta Math. Soc.}, 38:193--206, 1946.

\bibitem[Lah47]{MR0022566}
D.~B. Lahiri.
\newblock On {R}amanujan's function {$\tau(n)$} and the divisor function
  {$\sigma\sb k(n)$}. {II}.
\newblock {\em Bull. Calcutta Math. Soc.}, 39:33--52, 1947.

\bibitem[LR05]{GT0509205}
Samuel Leli\`evre and Emmanuel Royer.
\newblock Orbitwise countings in $\mathcal{H}(2)$ and quasimodular forms, 2005.
\newblock arXiv math.GT/0509205.

\bibitem[LW05]{LeWi05}
M.~Lemire and K.S. Williams.
\newblock Evaluation of two convolution sums involving the sum of divisor
  functions.
\newblock {\em Bull. Aust. Math. Soc.}, 73:107--115, 2005.

\bibitem[Mel98]{MR1628855}
Giuseppe Melfi.
\newblock On some modular identities.
\newblock In {\em Number theory (Eger, 1996)}, pages 371--382. de Gruyter,
  Berlin, 1998.

\bibitem[MR05]{MR2186573}
Fran{\c{c}}ois Martin and Emmanuel Royer.
\newblock Formes modulaires et p\'eriodes.
\newblock In {\em Formes modulaires et transcendance}, volume~12 of {\em
  S\'emin. Congr.}, pages 1--117. Soc. Math. France, Paris, 2005.

\bibitem[Ram16]{Ram16}
K.S. Ramanujan.
\newblock On certain arithmetical functions.
\newblock {\em Trans. Cambridge Philos. Soc.}, 22:159--184, 1916.

\bibitem[Ran56]{MR0082563}
R.~A. Rankin.
\newblock The construction of automorphic forms from the derivatives of a given
  form.
\newblock {\em J. Indian Math. Soc. (N.S.)}, 20:103--116, 1956.

\bibitem[Ser77]{MR0498338}
Jean-Pierre Serre.
\newblock {\em Cours d'arithm\'etique}.
\newblock Presses Universitaires de France, Paris, 1977.
\newblock Deuxi\`eme \'edition revue et corrig\'ee, Le Math\'ematicien, No. 2.

\bibitem[Shi72]{MR0314801}
Goro Shimura.
\newblock Class fields over real quadratic fields and {H}ecke operators.
\newblock {\em Ann. of Math. (2)}, 95:130--190, 1972.

\bibitem[Shi94]{MR1291394}
Goro Shimura.
\newblock {\em Introduction to the arithmetic theory of automorphic functions},
  volume~11 of {\em Publications of the Mathematical Society of Japan}.
\newblock Princeton University Press, Princeton, NJ, 1994.
\newblock Reprint of the 1971 original, Kan\^o Memorial Lectures, 1.

\bibitem[Ste04]{Ste04}
William Stein.
\newblock {Algorithms For Computing With Modular Forms}, 2004.
\newblock Preprint. Notes of a graduate course given at Harvard University,
  (Fall 2004).

\bibitem[Wil]{Wi}
K.S. Williams.
\newblock The convolution sum $\sum_{m<n/8}\sigma(m)\sigma(n-8m)$.
\newblock {\em Pacific J. Math.}
\newblock To appear.

\bibitem[Wil05]{MR2173379}
Kenneth~S. Williams.
\newblock The convolution sum {$\sum\sb {m<n/9}\sigma(m)\sigma(n-9m)$}.
\newblock {\em Int. J. Number Theory}, 1(2):193--205, 2005.

\bibitem[Zag92]{MR1221103}
Don Zagier.
\newblock Introduction to modular forms.
\newblock In {\em From number theory to physics (Les Houches, 1989)}, pages
  238--291. Springer, Berlin, 1992.

\end{thebibliography}
\end{document}